\documentclass[12pt,twoside]{article}

\usepackage{amssymb}
\usepackage{amsthm}
\usepackage{graphicx,color}
\usepackage{pst-all}
\usepackage{multido}
\def\R{{\mathbb R}}
\def\N{{\mathbb N}}
\def\C{{\mathbb C}}
\def\Z{{\mathbb Z}}

\newtheorem{thm}{Theora}[section]
\newtheorem{theo}[thm]{Theorem}

\newtheorem{lem}[thm]{Lemma}
\newtheorem{prop}[thm]{Proposition}
\newtheorem{Def}[thm]{Definition}

\newcommand\ds{\displaystyle\sum}
\newcommand\dso{\displaystyle \bigotimes}
\newcommand\dprod{\displaystyle\prod}
\newcommand\dmax{\displaystyle\max}
\newcommand\bi{\displaystyle\bigcup}

\def\prend{$~~\mbox{\hfil\vrule height6pt width5pt depth-1pt}$ }
\begin{document}
\pagestyle{myheadings} \markboth{ H. Gottschalk, B. Smii }{How to
determine the law of the noise driving a SPDE }
\thispagestyle{empty}

%\title{The Feynman graph representation of general convolution semigroups and its applications}
\title{How to
determine the law of the noise driving a SPDE }
 \author{Hanno Gottschalk${}^\sharp$ and  Boubaker Smii${}^\sharp {}^\flat$  }
\maketitle {\small

 \noindent ${}^\flat$: D\'epartement des Math\'ematiques, Universit\'e de Tunis El Manar

 \noindent  ${}^\sharp$: Institut f\"ur angewandte Mathematik, Rheinische Fridrich-Wilhelms-Universit\"at Bonn

\vspace{0.5cm} \begin{center} {\bf Abstract } \end{center} \small{We
consider a stochastic partial differential equation (SPDE) on a lattice
$$
\partial_t X=(\Delta-m^2)X-\lambda X^p+\eta$$
where $\eta$ is a space-time L\'evy noise. A perturbative (in the sense of formal power series)
strong solution is given by a tree expansion, whereas the correlation functions of the
solution are given by a perturbative expansion with coefficients that are represented as sums over a certain class of graphs, called Parisi-Wu graphs.
The perturbative expansion of the truncated (connected) correlation functions is obtained via a Linked Cluster
Theorem as a sums over connected graphs only. The moments of the stationary solution can be calculated as well.
In all these solutions the cumulants of the single site distribution of the noise enter as multiplicative constants. To determine them, e.g. by
comparison with a empirical correlation function, one can fit these constants (e.g. by the methods of least squares)  and
thereby one (approximately) determines law of the noise.\\}
%\tableofcontents

\vspace{1cm}

%$$(2\pi)^{-d}\int_{\Pi_\delta^d}e^{i{\bf k}\cdot ({\bf x}-{\bf y})}d{\bf
%k}=\delta^{-d}\delta_{{\bf x},{\bf y}}$$

%ou $\delta_{{\bf x},{\bf y}}$ est le symbol de Kronecker.

\newcommand{\GA}{\mbox{
\psset{xunit=1cm,yunit=1cm,runit=1cm,shortput=tab}
\begin{pspicture}(0,0)(.62,.3)
\dotnode[dotstyle=x,dotscale=2](0,0){D1}
\dotnode[dotstyle=otimes,dotscale=1.5](0.5,0){F1}\ncline{F1}{D1}
\end{pspicture}
}}
\newcommand{\GB}{\mbox{
\psset{xunit=1cm,yunit=1cm,runit=1cm,shortput=tab}
\begin{pspicture}(0,0)(.83,.3)
\dotnode[dotstyle=x,dotscale=2](0.,0){D1}
\dotnode[dotstyle=pentagon,dotscale=1.5](0.5,0){F1}\ncline{F1}{D1}
\end{pspicture}
}}
\newcommand{\GC}{\mbox{
\psset{xunit=1cm,yunit=1cm,runit=1cm,shortput=tab}
\begin{pspicture}(0,0)(.83,.3)
\ncarc[arcangle=45]{B}{C} \cnode*(0.5,0){2pt}{A}
\dotnode[dotstyle=x,dotscale=2](0.,0){D}\ncline{A}{D}
\dotnode[dotstyle=pentagon,dotscale=1.5](0.8,0.3){F}
\dotnode[dotstyle=pentagon,dotscale=1.5](1,0){F_1}
\dotnode[dotstyle=pentagon,dotscale=1.5](0.8,-0.3){F_2}
\ncline{A}{F} \ncline{A}{F_1} \ncline{A}{F_2}
\end{pspicture}
}}
\newcommand{\GD}{\mbox{
\psset{xunit=1cm,yunit=1cm,runit=1cm,shortput=tab}
\begin{pspicture}(0,0)(.62,.25)
\cnode*(0.5,0){2pt}{B}
\dotnode[dotstyle=otimes,dotscale=1.5](0.8,0.3){G}
\dotnode[dotstyle=otimes,dotscale=1.5](1,0){G_1}
\dotnode[dotstyle=otimes,dotscale=1.5](0.8,-0.3){G_2}
\dotnode[dotstyle=x,dotscale=2](0.,0){H}\ncline{B}{H}\ncline{B}{G}\ncline{B}{G_1}\ncline{B}{G_2}
\ncline{D}{A}
\end{pspicture}
}}
\newcommand{\GE}{\mbox{
\psset{xunit=1cm,yunit=1cm,runit=1cm,shortput=tab}
\begin{pspicture}(0,0)(1.23,.3)
\cnode*(0.5,0){2pt}{B_1}
\dotnode[dotstyle=pentagon,dotscale=1.5](0.8,0.3){G'}
\dotnode[dotstyle=otimes,dotscale=1.5](1,0){G''}
\dotnode[dotstyle=pentagon,dotscale=1.5](0.8,-0.3){G'''}
\dotnode[dotstyle=x,dotscale=2](0.,0){H_1}\ncline{B_1}{H_1}\ncline{B_1}{G'}\ncline{B_1}{G''}\ncline{B_1}{G'''}
\end{pspicture}
}}
\newcommand{\GF}{\mbox{
\psset{xunit=1cm,yunit=1cm,runit=1cm,shortput=tab}
\begin{pspicture}(0,-.05)(.83,.3)
\cnode*(0.5,0){2pt}{B_1}
\dotnode[dotstyle=otimes,dotscale=1.5](0.8,0.3){G'}
\dotnode[dotstyle=pentagon,dotscale=1.5](1,0){G''}
\dotnode[dotstyle=otimes,dotscale=1.5](0.8,-0.3){G'''}
\dotnode[dotstyle=x,dotscale=2](0.,0){H_1}\ncline{B_1}{H_1}\ncline{B_1}{G'}\ncline{B_1}{G''}\ncline{B_1}{G'''}
\end{pspicture}
}}
\newcommand{\GH}{\mbox{
\psset{xunit=1cm,yunit=1cm,runit=1cm,shortput=tab}
\begin{pspicture}(0,-.09)(.62,.25)
$\hspace*{3mm}\Big)+o(\lambda^2)$
\end{pspicture}
}}

\section{Introduction and overview}

In this work we give a sufficiently general perturbative solution to the Stochastic Partial Differential Equation
(SPDE)
 \begin{equation} \label{1.1eqa}
\left\{ \begin{array}{ll}{{\partial X}(t,\,{\bf x}) \over {\partial
t}}= \Delta X(t,\,{\bf x}) -m^2 X(t,\,{\bf x})- \lambda X^p(t,\,{\bf
x})+\eta(t,\,{\bf x})&, \lambda
> 0,\,m>0, p\in \N
\\ X(0,\,{\bf x})=f({\bf x}),\,\,\,(t,\,{\bf x}) \in ]0,\,\infty[ \times L_{\delta}
\end{array}
\right.
\end{equation}
where $\eta(t,{\bf x})$ is a general space-time noise of L\'evy type with moments of arbitrary order such that such that one can extract information on the
distribution of $\eta$ by comparison with empirical data. To avoid so-called ultra-violet (UV) divergences, we consider this equation on a
spacial lattice $L_{\delta}=\{\delta z\,, z \in \Z^d\}$ with lattice spacing $\delta>0$.

Nonlinear SPDEs driven by non Gaussian noise have been discussed
both from the physical \cite[Section 4.2 and Refs.]{HHZ} and the
mathematical side, see e.g. \cite{AWZ,AW}. But so far there seems to
be no recipe to find out which kind of noise -- given the general
structure of the equation, including, say, the numerical values of
all coefficients the right hand side -- explains a set of empirical
data best. It is the aim of this article to give such a recipe by
calculating the correlation (moment) functions of the solution $X$
perturbatively. To do so, we generalize the classical Feynman graph
approach of Parisi and Wu \cite{PW} on the Gaussian stochastic
quantization equation to the general L\'evy case, using a
generalized class of graphs, henceforth called Parisi-Wu graphs.
Statistically relevant quantities as expectation, variance and
higher order cumulants of the noise enter in the expansion as simple
constants -- like additional coupling constants --  and can
therefore be extracted rather easily by comparison with the data.

The restriction to equation (\ref{1.1eqa}) is not essential at all
for the combinatorial graph calculus that we develop. The following
generalizations are rather obvious: i) the nonlinearity $X^p$ can be
replaced by a polynomial, e.g. the Mexican hat potential $2X(X^2-a)$
in equations of reaction-diffusion type; ii) nonlinearities include
derivatives, like $|\nabla X|^2$ in the KPZ-equation \cite{BS,HHZ};
iii) in addition the linear operator content on the right hand side
is of higher order and the noise is of derivative type as in the
non-linear description of nuclear beam epitaxy, where the right hand
side of (\ref{1.1eqa}) is $-K\Delta^2X-\lambda \Delta |\nabla
X|^2+\Delta\eta$ \cite{BS,HHZ}; iv) the left hand side is of second
order $\partial_t^2X$, as in models of thermalization and symmetry
breaking in the early universe \cite{BF}. All these changes can be
incorporated in our formalism by adaptation of the so-called Feynman
rules without changing the combinatorial core of this paper. Hence,
it is mostly for notational convenience that here we only consider
the SPDE (\ref{1.1eqa}).

A perturbative expansion in the sense of formal power series without
any control on the convergence of the series, as done in this work,
can certainly be criticized from a mathematical point of view. In
fact, even the existence of a solution to (\ref{1.1eqa}) for the
general L\'evy case to our best knowledge has not been studied
though the restrictions on the nonlinearities in the continuum
\cite{AW} that exclude polynomial interactions are related to the
above mentioned UV-singularities and most likely can be dropped in
the lattice regularized case. In this work we happily follow the
tradition of perturbative quantum field theory (QFT) that, in spite
of this, such expansions  deliver (after suitable resummation
\cite[Chapters 42.5--42.7]{ZJ}) highly precise numerics.

The paper consists of seven sections.  In the next section we
give a solution of (\ref{1.1eqa}) in the sense of formal power
series in the parameter $\lambda$
\begin{equation}\label{b1}X(t,\,{\bf x})=\ds_{j=0}^\infty (-\lambda)^j X_j(t,{\bf x}).\end{equation}
Here $X_0$ is the solution of the linear equation ($\lambda=0$) and the $X_j$, under suitable conditions on the stochastic driving force $\eta$, are determined recursively.

In section 3 we represent this solution as a sum over rooted trees with two types of leaves, cf. Appendix A for the graph theoretic notions used in this article.
One type of leave is standing for the noise
$\eta$ and the other for the initial condition, cf. the following first order example for $p=3$
\begin{eqnarray}
\label{1.Eeqa} X(t,\,{\bf
x})&=&\GA+\GB+\lambda\,\,\Big[\GC\,\,\,+\,\GD\,\,\,\,\,\,\,+\,\,\,3\,\,\,\GE\,+\,\,3\,\,\GF\,\,\,\Big]+\,
o(\lambda^{2})\hspace*{6mm}
\end{eqnarray}
In section 4 we then develop the graphical calculus for the
perturbative evaluation of the correlation functions, generalizing
what in physics is called the Parisi-Wu method from Gaussian to
L\'evy noise. While only the well-known Feynman graphs occur in the
Gaussian case, we have to pass on to a more general class of graphs
in order to deal also with L\'evy type noise. We also prove a Linked
Cluster Theorem that shows that the truncated correlation functions
(cumulants of the process $X$) are given by the sum over connected
Parisi-Wu graphs, only. So far no detailed assumptions on the
stochastic driving force were needed.

In section 5 we show how the rules for the analytic evaluation of
graphs simplifies if the noise $\eta$ is to be taken as a white
noise of L\'evy type. We also show that a L\'evy noise $\eta$
fulfills the requirement for the existence of the perturbative
solution as specified in section 2.

In section 6 we study  the limit where $t$ goes to infinity and give the perturbative expansion of equilibrium correlation functions.

In section 7 we briefly study the question how to extract statistical information
on the noise from an empirical correlation function, that e.g. could be obtained by X-ray analysis or sampling
methods. To illustrate this point we solve a least square minimization problem to first order perturbation theory and give some indications, how to perform
 higher order calculations.

\section{ Formal solution of the SPDE}

In this section we give a formal solution of stochastic partial
differential equation (\ref{1.1eqa}).\\ We first fix some notation.
Let $\Gamma=\R \times L_\delta\,=\{x=(t,\,{\bf x})\,\,\,;\,\,\,t
\,\in \R\,\,,\,\,{\bf x}\,\in \,L_\delta \}$ and $\Pi_\delta^d \cong
[0,\,  \,{{2\,\pi} \over \delta}]^d$ where the opposite edges are
being identified. Let $B=L_\delta$ or $B=\Gamma$. By $S(B)$ we
denote the Schwartz Space of all rapidly decreasing functions on $B$
endowed with the Schwartz topology, its topological dual is the
Space of tempered distribution noted by $S'(B)$. We denote by
$\langle.\,,\,.\rangle$ the dual pairing between
$S(B)$ and $S'(B).$\\
Now we introduce the following notation :
For $f \,\in\, S(\Gamma)$, $\int f(x)\,d x=\sum_{x \in L_{\delta}} \delta^d
\int_{\R}f(t,\,{\bf x})dt$, $x=(t,\,{\bf x})$. For $f
\,\in\,S(L_\delta)$, $\int f({\bf x}) d{\bf x}=\sum_{x \in L_{\delta}} \delta^d f({\bf
x})$. Furthermore, for $A\subset \Gamma$ ($A\subset L_\delta$), $\int_A f(x) dx=\int 1_A(x)f(x) dx$ ($\int_A f({\bf x}) d{\bf x}=\int 1_A({\bf x})f({\bf x}) d{\bf x}$) where $1_A$ is the characteristic function of the set $A$ and $f\in S(\Gamma)$ ($f\in S(L_\delta)$). \\\\Let ${\cal{F}}:
S(\Gamma)\,\longrightarrow\,\R \times
\Pi_\delta^d$ be the Fourier transform on $S(\Gamma)$ \\
i.e.\begin{equation}{\cal F}(f)(E,\,{\bf p})=\int e^{itE}e^{i{\bf
{p.x}}}f(x) dx\,\,,\,\,f \,\,\in \,\,S(\Gamma).\end{equation}The
inverse Fourier transform of, ${\cal F}$, is given by
\begin{equation}f(t,\,{\bf x})={1 \over
{(2\pi)^{d+1}}}\int_{\R}\int_{\Pi_\delta^d} e^{-iEt}e^{-i{\bf p\cdot
x}}{\cal F}(f)(E,\,{\bf p})dE\,d{\bf p}.\end{equation} and let ${
\tilde{\cal F}}:S(L_\delta)\,\longrightarrow\, \Pi_\delta^d$ be the
Fourier transform on $S(L_\delta)$\\i.e :\begin{equation}{
\tilde{\cal F}}(f)({\bf p})=\int e^{i{\bf {p.x}}}f({\bf x})d{\bf x}
\,\,,\,\,f \,\in\,\,S(L_\delta).\end{equation} The inverse Fourier
transform of, $\tilde{{\cal F}}$ is given by :
\begin{equation}f({\bf x})={1 \over {(2\pi)^d}}\int_{\Pi_\delta^d}
e^{-i{\bf p}.{\bf x}}\tilde{{\cal F}}(f)({\bf p})d{\bf
p}.\end{equation}By definition the action, of the Lattice Laplacian
$\Delta$ on a test function $f\,\in\,S(L_\delta)$ is as follows :
\begin{equation}\Delta f({\bf x})=\delta^{-2}[-2d\, f({\bf
x})+\ds_{\mid {\bf x}-{\bf y}\mid =\delta}f({\bf y})].\end{equation}
Then
\begin{equation} \tilde{{\cal F}}(\Delta f)({\bf
p})=\delta^{-2}\Big[-2d+2\ds_{j=1}^d
\cos(\delta {\bf p}_j)\Big]\tilde{{\cal F}}(f)({\bf p}).
\end{equation}
Where $e_j=(0,...,1,0,...)$ is the canonical basis of $\R^d$ and
${\bf p}_j={\bf p}{\cdot}e_j.$
\\We define
\begin{equation}\label{b2}\mu^2_{\delta,\,m}({\bf p})=2\delta^{-2}\Big(d-\ds_{j=1}^d
\cos(\delta {\bf p_j}) \Big)\,+{m^2 }.%=\nu^2_{\delta,\,m}({\bf
%p})+{m^2 \over 2}.
\end{equation}
We introduce two convolutions product $"\ast"$ and $"\star"$
respectively by 
\begin{equation}f\ast g( x)=\int f( x- y)g( y)\, dy,\,\,\,f,\,g\,\, \in
\,\,S(\Gamma).\end{equation}
and
\begin{equation} f\star g({\bf x})= \int f({\bf x}-{\bf y}) g({\bf
y})\,d{\bf y}\,\,,\,\,f,\,g\,\in \,S(L_\delta).\end{equation}
Let $G(t,\,{\bf x})$ be the
Green function which satisfies :
\begin{equation} \label{1.11eqa}
\left\{ \begin{array}{ll}{{\partial G}(t,\,{\bf x}) \over {\partial
t}}= \Delta G(t,\,{\bf x}) -m^2 G(t,\,{\bf x})+ \delta( x)&,
\\ G(t,\,{\bf x})=0,\,\,\, \,\,\,t<0\,\,\,\,\,\,\,\,\,
\end{array}
\right.
\end{equation}
Here $\delta(x)$ is the Dirac distribution defined by $\delta(x)=\delta(t)\delta({\bf x})=\delta(t)\delta^{-d}\delta_{0,\,{\bf x}}$ where $\delta(t)$ is the Dirac distribution on $\R$ and $\delta_{{\bf x},{\bf y}}=\prod_{l=1}^d\delta_{\delta^{-1}x_l,\delta^{-1}y_l}$ with $\delta_{i,j}$ the Kronecker symbol.\\
Applying to the Fourier transform ,${\cal F}$, it is obvious that
\begin{equation}\label{equa1}{\cal{F}} (G)(E,\,{\bf p})= {1 \over {-iE +
\mu^2_{\delta,\,m}({\bf p})}}\,\,\,.\end{equation} Let
$\tilde{G}_t({\bf x})$ be the Green function which satisfies :
\begin{equation} \label{1.111eqa}
\left\{ \begin{array}{ll}{{\partial\tilde{G}_t}(\,{\bf x}) \over
{\partial t}}= \Delta \tilde{G}_t(\,{\bf x}) -m^2 \tilde{G}_t(\,{\bf
x})&,
\\ \tilde{G}_0(\,{\bf x})=\delta({\bf x}).\,\,\,\,\,\,\,\,\,\,\,\,
\end{array}
\right.
\end{equation}
Then by application of the Fourier transform $\tilde{{\cal F}}$ we
get 
\begin{equation} \label{1.12eqa} \tilde{\cal
F}(\tilde{G}_t)({\bf p})=e^{-t\mu_{\delta,\,m}^2({\bf
p})}\,\,\,.\end{equation} By the use of the residue Theorem, one can
see that $G$ and $\tilde{G}_t$ are related by the following equation
\begin{equation}
G(t,\,{\bf x})=\theta(t)\tilde{G}_t(\,{\bf x}),\end{equation} where
$\theta(t)=1$ if $t>0$ and $\theta(t)=0$ else. Let
$\Delta_{\Pi^d_\delta}$ be the Laplacian with cyclic boundary
condition on $[0,\frac{2\pi}{\delta}]\cong\Pi_\delta^d$. We get
\begin{lem}\label{l.1}$\forall \,\,N \in \N,\,\,\exists\,\,K_N $ such that :$$ \mid G(t,\,{\bf x})\mid \leq\,K_N\,\,
{e^{-{tm^2 \over 2}} \over {(1+\mid{\bf
x}\mid^2)^N}},\,\,\,\,(t,\,{\bf x})\,\in\,\,
]0,\,\infty[\,\times\,L_\delta.$$
\end{lem}
\noindent {\bf Proof.} It suffices to prove \begin{equation}
\sup_{{\bf x}} \mid G(t,\,{\bf x})\,(1+\mid{\bf x}\mid^2)^N \,\mid
\,<\,K_N\,e^{-{tm^2 \over 2}}.\end{equation}
\begin{eqnarray}
G(t,\,{\bf x})\,(1+\mid{\bf x}\mid^2)^N \,&=&{1 \over
{(2\,\pi)^d}}(1+\mid{\bf x}\mid^2)^N\theta(t)\,\int_{\Pi ^d_\delta}e^{i{\bf
p}\cdot{\bf x}} e^{-t\mu_{\delta,\,m}^2({\bf p})}\,d{\bf
p}\nonumber\\&=&{1 \over {(2\,\pi)^d}}\theta(t)\int_{\Pi
^d_\delta}\left[(1-\Delta_{\Pi^d_\delta})^N\,e^{i{\bf p}\cdot{\bf x}}\right]
e^{-t\mu_{\delta,\,m}^2({\bf p})}\,d{\bf
p}\nonumber\\&=&{1 \over {(2\,\pi)^d}}\theta(t)\int_{\Pi ^d_\delta}\,e^{i{\bf
p}{\bf x}}
\left[(1-\Delta_{\Pi^d_\delta})^N\,e^{-t\mu_{\delta,\,m}^2({\bf
p})}\right]\,d{\bf p}.
\end{eqnarray}
Then \begin{equation}\sup_{\,x\in L_\delta} \left| G(t,\,{\bf
x})\,(1+\mid{\bf x}\mid^2)^N \,\right|\,\leq\,e^{-{tm^2 \over
2}}\,\left({2\pi\over\delta}\right)^{d}\,\,
 \sup_{t\in (0,\infty),\,{\bf p}\in \Pi^d_\delta} \left| (1-\Delta_{\Pi^d_\delta})^N \,
\,e^{-t\mu^2_{\delta,\,{m\over \sqrt{2}}}({\bf
p})}\right|.\end{equation} Now take
\begin{equation}K_N=\left({2\pi\over\delta}\right)^{d}\,\,\sup_{t,\,{\bf
p}} \left| {(1-\Delta_{\Pi^d_\delta})^N\,e^{-t\mu^2_{\delta,\,{m
\over \sqrt{2}}}({\bf p})}}\right|.\end{equation}\\{\hspace*{15cm
\prend}}
\begin{Def}
{\rm We define a space of all measurable and polynomially bounded
functions by :\begin{equation}{\cal M}_b=\left\{f:\,
\Gamma\,\longrightarrow\,\R\ \mbox{\rm measurable};\,\forall\,k\,\in
\N,\,\,\exists\,\,N\,\;\in \,\N\,\,:\,\,\int\,{|f|^k(x)\over
\,(1+\mid x\mid^2)^N}\,\,dx<\,\infty\right\}\end{equation}
}\end{Def}
\begin{lem}
\label{1.2alem}
${\cal M}_b$ is an Algebra under multiplication.
\end{lem}
\noindent {\bf Proof.} For $f,\,g\,\,\,\in \,{\cal M}_b$, $fg$ and
$f+ag$, $a\in\R$, is a measurable function. The vector space
structure of ${\cal M}_b$ easily follows from the triangular
inequality.

Furthermore, for $k\in\N$ fixed, $\exists\,\,N\,\in \,\N$ such that
\begin{equation}\int\,{f^{2k}(x)\over \,(1+\mid
x\mid^2)^{N}}\,\,dx<\,\infty~~~\mbox{and}~~~\int\,{g^{2k}(x)\over
\,(1+\mid x\mid^2)^{N}}\,\,dx<\,\infty\end{equation} Due to the
Cauchy Schwartz inequality,
\begin{equation}
\int\,{|f g|^k(x)\over \,(1+\mid x\mid^2)^N}\,\,dx
\leq\sqrt{\int\,{f^{2k}(x)\over \,(1+\mid
x\mid^2)^{N}}dx\,\int\,{g^{2k}(x)\over \,(1+\mid
x\mid^2)^{N}}dx}~~<\infty. \end{equation}
\prend
\begin{lem}\label{1.3alem}
Let $f\,\in \,{\cal M}_b$ then $G\ast\,f\,\,\in\,{\cal M}_b.$
\end{lem}
\noindent{\bf Proof.} Let $f\in{\cal M}_b$. We have due to Lemma \ref{1.1lem}
\begin{eqnarray}
\int{|G\ast f|^k(x)\over\,(1+\mid{
x}\mid^2)^{2kN}}\,\,dx &\leq& C(N)\,\int\,{{\int\,{|f|(y_1)\over(1+\mid{
x- y_1}\mid^2)^{N}
}\,dy_1\cdots\int{|f|(y_k)\over(1+\mid{ x-
y_k}\mid^2)^{N}}dy_k}\over(1+\mid{ x}\mid^2)^{2kN}}\,\,dx
\nonumber\\&=&C(N)\int\,{\dprod_{l=1}^k|f|(y_l)\over{\dprod_{l=1}^k\left[(1+\mid x- y_l\mid^2)^{N}(1+\mid{
x}\mid^2)^{2N}\right]}}\,\,dx\,dy_1\cdot\cdot\cdot
dy_k\nonumber\\&\leq&\,4^kC(N)\,\left(\,{\int\,{|f|(y)\over(1+\mid{
y}\mid^2)^{N}}}\right)^k\,\int{dx\,\over {(1+\mid
x\mid^2)^{N k}}}\,.
\end{eqnarray}
Here the right hand side is finite for $N$ sufficiently large. The last inequality is due to the fact that 
\begin{equation}
4(1+\mid x-y\mid^2)(1+\mid x\mid^2)^{2}
\geq(1+\mid x\mid^2)(1+\mid y\mid^2)\,\,\forall\,x,\,y\,\in
\Gamma.
\end{equation}
{\hspace*{15cm \prend}}
\begin{prop}\label{1.1prop}
Suppose that $G*\eta\in{\cal M}_b$ a.s.\footnote{For L\'evy noise with all moments this will be verified explicitly, cf. Proposition \ref{PROP}}.  Let $X(t,\,{\bf x})=\sum_{j=0}^\infty (-\lambda)^j X_j(t,\,{\bf x})$
be an expansion of $X$ in the sense of formal power series in the parameter
$\lambda$. The perturbative solution of the stochastic differential equation
(\ref{1.1eqa}) is given by \begin{equation} \label{1.2eqa}
\left\{ \begin{array}{ll} X_0(t,\,{\bf x})=G \ast \eta(t,\,{\bf x})+\tilde{G}_t\star f({\bf x})\\
X_j(t,\, {\bf x})=G\ast \ds_{B^j_p({n_0,\,n_1,...})} {p! \over
{n_0!n_1!...}} \dprod_{i\geq 0} X^{n_i}_i(t,\,{\bf x}),\,\,\,j\geq
1.
\end{array}
\right.
\end{equation}
with $X_j\in{\cal M}_b$ a.s.. Here $B^j_p(n_0,\,n_1,\,...)=\{n_0,\,n_1,\,...\geq \,0 \,\,\mid\,\,
\, n_0+n_1+...=p\,;\,\,\ds_{i\,\geq\,0}i\,n_i\,=j-1\,\}$.
\end{prop}
\noindent {\bf Proof.} We have the integrated form of the SPDE
\begin{equation}X(t,\,{\bf x}) =  - \lambda G*(X^p)(t,\,{\bf x})+G*\eta(t,\,{\bf
x})+\tilde G_t\star f({\bf x})\end{equation} It is obvious that the
solution of this equation for the zero-th order in $\lambda$ is
given by :
\begin{equation} \label{1.11eqa}X_0(t,\,{\bf x})=G\ast \eta(t,\,{\bf x})+\tilde{G}_t\star
f({\bf x}),\end{equation} where $G$ and $\tilde{G}$ are the Green
functions defined earlier.
Now we determine the solution of higher order in $\lambda$. We have
\begin{eqnarray} \ds_{j=0}^\infty (-\lambda)^j X_j(t,\,{\bf x})&= &- \lambda
\, G*\left(\ds_{j=0}^\infty (-\lambda)^j X_j\right)^p(t,\,{\bf
x})+G*\eta(t,\,{\bf x})+\tilde G_t\star f(\,{\bf
x})\nonumber\\&=&-\lambda \,\,\,G*\left(
\ds_{{n_0,\,n_1,...\geq 0} \atop {n_0+n_1+...=p}} {p! \over
{n_0!n_1!...}} \dprod_{i\geq 0} \lambda^{in_i} X^{n_i}_i (t,\,{\bf
x})\right)\nonumber\\&&~~~~~~~~~~~~~~~~~~~~~~~~~~~~~~~~~~~+G*\eta(t,\,{\bf x})+\tilde G_t\star f(\,{\bf
x}).\nonumber\end{eqnarray}
By comparison of coefficients we get \\
\begin{equation}\label{1.000eqa}X_j(t,\,{\bf x})=G\ast \ds_{{n_0,\,n_1,...\geq 0} \atop{
{n_0+n_1+...=p}\atop{0\times n_0+1\times n_1+...={j-1}}}} {p! \over
{n_0!n_1!...}} \dprod_{i\geq 0}
X^{n_i}_i(t,\,{\bf{x}})\,\,\,\,;\,\,\,\,j\geq\,1.\end{equation}
$X_0\,\in \,{\cal M}_b$ a.s. by our assumptions. Using Lemma \ref{1.2alem} and Lemma \ref{1.3alem}
it follows by induction that $X_j\,\,\in \,\,{\cal M}_b$ a.s. $\forall\,\,j\,\in
\N$.\prend

\section{Tree expansion of the perturbative  solution}
\hspace*{5mm}In this section we first recall some notions of graph
theory. We then recall a fundamental result which represents the
solution -- in the sense of a formal power series -- of the stochastic differential equation
(\ref{1.1eqa})  by a sum over all rooted trees.
Some of the graph theoretic standard notions, that are used here, are collected in Appendix A and also \cite{JPS,OYO}.\\
A tree $T$ is a connected graph  without cycle. We also consider graphs with different types of vertices. A rooted tree with root
$x \,\in\, \Gamma$, denoted by $\times$, and two types of leaves, the
leaves of type one are noted by $\otimes$ , and those of type two
by ~\begin{pspicture}(0,0)(0,0)\psdot[dotstyle=pentagon,dotsize=7pt](0,.1)\end{pspicture}\hspace{.4cm}.\\
We define a partial order $"\leq"$ on the set of the vertices of
$T$, $V(T)$ by $v \leq w$ if every walk connecting $w$ and
$x$ passes through $v$. \\
We note ${\cal A}(i)$ the set of all rooted trees $T$ with root $x$
and two types of leaves, which has $i$ inner vertices with $p+1$
legs. See e.g. (\ref{1.Eeqa}) for the trees from ${\cal A}(0)$ and $\cal{A}(1)$ for $p=3$.

\begin{Def}{\label{1.1def}}
{\rm For $T \in {\cal A}(i)$, the random variable ${\cal
B}(T,\,\eta,\,x),$ is defined as follows:\\1)-Assign $x \in \Gamma$ to the root of the tree $T$.\\
 \hspace*{2mm}- Assign  values $y_1,...,y_i \in \Gamma$ to the inner vertices.\\
 \hspace*{2mm} -Assign  values $z_1,...,z_k \in \Gamma\,,\, $ to the
leaves of type one and assign a values ${\bf z}'_1,...,{\bf z}'_l
\in L_{\delta}$ to the leaves of type two.\\
2)- For every edge with two end points,$e=\{v,\, w\},\,\,$ assign a
value $G(e)=G(v-w),\,\,\,(v\leq w)$ to this edge. $G$ is the
Green function defined in (\ref{equa1}). \\
3)- For the $j$-th leaf multiply with $\eta(t_j,\,{\bf z}_j)$ if
this leaf is of type one and multiply
with $f({\bf z}_j)$ if this leaf is of type two.\\
4)- For the $l-th$ inner vertex multiply with the coefficient
${p!\over {n'_0!(n_0-n'_0)!n_1!...}},$ where $n_i,\, i \in \N,$ is
the number of rooted subtrees connected to the $l-th$ inner vertex
with $i$ inner vertices and $n'_0$ is the number
of the rooted subtrees connected to this vertex with zero inner vertices and a leaf of type one.\\
5)- Integrate with respect to the Lebesgue measure
$dy_1\cdot\cdot\cdot dy_p\,dz_1\cdot\cdot\cdot dz_k\,d{\bf
z}'_1\cdot\cdot\cdot d{\bf z}'_l$.
}\end{Def}

Let $T \in {\cal A}(j)$, $j \in \N$, and
$C(j)=\{T_{0,1},...,T_{0,\,n_0}T_{1,1},...;\,T_{i,1},...,T_{i,\,n_i}\in
{\cal A}(i)\,,\,\sum_{i\geq 0}in_i=j-1,\,\,n_0+n_1+...=p\}$. We
construct a one to one correspondence between ${\cal A}(j)$ and
$C(j)$ that is given in the following way:\\
We first consider the case when $T\,\in\,{\cal A}(j) $ is
given.
Let $y_1$ be the first inner vertex of $T$ and we denote
$e_1=\{x,\,y_1\}$ the first edge of $T.$ Cut the tree $T$ at
$y_1$. Let $n_i,\,i\,\geq\,0$ be the number of the rooted subtrees
with root $y_1$ and $i$ inner vertices denoted by
$T_{i,1},...,T_{i,{n_i}}.$ Then $n_0+n_1+...+n_i+...=p $ and
$\sum_{i\geq 0}i\,n_i=j-1$. Conversely, let $x \in \Gamma$ and
$T_{i,1},...,T_{i,n_i} \in {\cal A}(i)$, be a collection of rooted
trees with root $y_1$ be given such that $n_0+n_1+...=p
,\,\sum_{i\geq 0}i\,n_i=j-1$. Attach this collection of rooted trees
to the edge $e=(x,\,y_1)$.
By definition we obtain a rooted tree $(T,\,x) \in {\cal A}(j)$. For the notions of cutting and attaching we again refer to Appendix A.\\
The second operation is obviously the inverse of the first one. We
have established the following result:
\begin{lem}\label{1.1lem}
Let $T \in {\cal A}(j), (j \in \N)$ and $T_{i,1},...,T_{i,\,n_i} \in
{\cal A}(i),\,i \in \N$, the application \begin{equation}{\cal
C}:{\cal A}(j)\longrightarrow C(j)\end{equation}is a one to one
correspondence and \begin{equation}{\cal A}:C(j)\longrightarrow
{\cal A}(j)\end{equation}is its inverse.
\end{lem}
\begin{Def}\label{def11}
{\rm Let $T\,\in \,{\cal A}(i)$ and $v\,\in\,V(T),$ we define
the multiplicity of the vertex $v$ as $M(v)={p!\over{n'_0!(n_0-n_0')!n_1!...}},$
where $n_i\,i\geq0,\,$ is the number of the rooted subtrees connected
to $v$ with $i$ inner vertices.\\The multiplicity of the rooted tree
$T$ is $M(T)=\dprod_{v \in V(T)}M(v).$
}\end{Def}
For $f\,\in \,S(L_\delta)$ we put $F_f(x)=\delta_0(t)f({\bf
x})\,\,\,,\,\,\,\,x=(t,\,{\bf x})\,\in\, L_\delta$
\begin{lem}\label{1.2lem}
Let $T \in {\cal A}(j)$ and $v_0 \in V(T)$ the vertex connected to
the root $x$, then the following result hold :\begin{equation}{\cal
B}(T,\,x,\,\eta)= M(v_0)G\ast[\dprod_{S \in{ \mathcal{C}(T)}}{\cal
B}(S,\,.\,,\,\eta) ](x). \end{equation}Where $\mathcal{C}(T)$ is the
cut (cf. Appendix A) of the  tree $T$.
\end{lem}\noindent{\bf Proof.}
\begin{eqnarray} \label{1.5eqa} {\cal
B}(T,\,x,\,\eta)&=&M(T)\int \dprod_{e \in E(T)} G(e) \dprod_{l_1 \in
L_1(T)} \eta(l_1)\dprod_{{ l} \in L_2(T)}F_f({l})
\dso_{v \in V(T)\setminus \{x\}}dv \nonumber\\&=& M(v_0)\int
G(x-v_0)\dprod_{S \in \mathcal{C}(T)}\Big[M(S) \int\dprod_{e \in
E(S)}G(e)\dprod_{l_1 \in L_1(S)}\eta(l_1)\nonumber\\&\times&\dprod_{
{ l} \in L_2(S)} F_f({l}) \dso_{v \in {V(S)\setminus {\{v_0\}}}} dv
\Big]dv_0 =M(v_0)G\ast[\dprod_{S \in
\mathcal{C}(T)}{\cal B}(S,\,.,\,\eta) ](x).
\end{eqnarray}
{\hspace*{15cm \prend}} \\Using Proposition {\ref{1.1prop}} and
Definition {\ref{1.1def}} we get the following theorem:
\begin{theo}\label{1thm1}
The solution of the stochastic differential equation $(1)$ in the
sense of formal power series is given  by a sum over all rooted
trees $T \in {\cal A}(j)$ that are evaluated according to the ruled
fixed in Definition \ref{1.1def}
i.e.\begin{equation}\label{a1eq}X_j(x,\,\eta)=\ds_{T \in {\cal
A}(j)}{\cal B}(T,\,x,\,\eta)\,\,\,\,\,;\,\,\,\,x\,\in
\,\Gamma\end{equation}
\end{theo}
\noindent{\bf Proof.} We prove the assertion by induction. For
$j=0$, we have $X_0(x,\,\eta)=G\ast\eta(x)+G\ast F_f(x)$, which just is the sum of the evaluation of the two trees in ${\cal A}(0)$ .\\ Suppose that ($\ref{a1eq}$) is true for all $i<j$. By
proposition $(2.2)$ we have \begin{equation}X_j(x,\,\eta)=G\ast
\ds_{B^j_p({n_0,\,n_1,...})} {p! \over {n_0!n_1!...}}\,\,
\dprod_{i\geq 0} X^{n_i}_i(x,\eta)\,\,  \end{equation} Then using
the induction hypothesis and Lemmas \ref{1.1lem} and
\ref{1.2lem} we get
\begin{eqnarray} \label{1.6eqa}
 X_j(x,\,\eta)&=&G\ast \ds_{B^j_p({n_0,\,n_1,...})} {p! \over {n_0!n_1!...}}
\,\,\dprod_{i\geq 0} \left[\ds_{T_i \in {\cal A}(i)}{\cal
B}(T_i,\,.,\,\eta)\right]^{n_i}(x)\nonumber \\&=& \ds_{B^j_p({n_0,\,n_1,...})} \sum_{T_{i,1},\ldots,T_{i,n_i}\in{\cal A}(i)\atop i=0,1,2,\ldots}{p!
\over {n_0!n_1!...}} \left(G \ast \dprod_{i\,\geq\,0}
\dprod_{l=1}^{n_i} {\cal B}(T_{i,\,l},\,.\,,\,\eta)\right)(x)\nonumber\\
&=&\ds_{B^j_p({n_0,\,n_1,...})} \sum_{T_{i,1},\ldots,T_{i,n_i}\in{\cal A}(i)\atop i=0,1,2,\ldots}{\cal B}({\cal A}(T_{0,1},\ldots,T_{0,n_0},T_{1,1},\ldots),x,\eta)\nonumber\\
&=& \sum_{T\in{\cal A}(j)}{\cal B}(T,x,\eta).
\end{eqnarray}
\hspace*{15cm} \prend
\section{A graphical calculus for the correlation functions}
In the following we define a class of graphs that we  call Parisi-Wu graphs. We then give a rule
which assigns a value to such a graph.  We prove a theorem which
represents the moments of the solution given in the previous
section in terms of a sum over all generalized Parisi-Wu graphs.
We finally  prove a linked cluster theorem for Parisi-Wu graphs.
\begin{Def}
{\rm 
A generalized Parisi-Wu graph of $m$-th order with $n$ exterior vertices
$x_1,...,x_n$ of fertility one, is a graph that contains $n$ rooted trees $T_1 \in {\cal
A}(j_1),...,T_n \in {\cal A}(j_n)$ as subgraphs, such that any leave of type one
is connected to exactly one vertex of a new type, called inner empty, of arbitrary fertility.
We denote such a graph by $\cal{G}$ and the set of all
generalized Parisi-Wu
graphs of $m$-th order and $n$ roots by $P(m,\,n)$.\\
If we denote  by $Q=Q({\cal G})$ the set of the empty vertices then
\begin{equation}V({\cal{G}})=\bi_{k=1}^n V(T_k)\cup Q
,\,\,\,E({\cal{G}})=\bi_{k=1}^nE(T_k)\cup E(T_k,\,Q)\end{equation}
where $E(T_k,\,Q)$ is a set of edges $e=\{v,\,q\},\,\,v \in
L_1(T_k),\,\,q \in Q$ such that $p(v)=2,\, \forall\, v \in L_1(T_k)$.\\
By definition, inner vertices are distinguishable and have non
distinguishable legs whereas empty vertices are non distinguishable
and have non distinguishable legs\footnote{The notion  of distinguishable and non distinguishable legs
 are explained in Appendix A.}.
}\end{Def}
The name generalized Parisi-Wu graph has been chosen to distinguish the graphs used here from Gaussian Parisi-Wu graphs implicitly defined in \cite{PW}, which, in a topological sense, coincide with the ordinary Feynman graphs.\\

\begin{prop}
For $j_1,...,j_n \in \N$, let $T_1 \in {\cal A}(j_1),...,T_n \in
{\cal A}(j_n)$ $n$  rooted trees with roots
$x_1,...,x_n\,\,\in\,\Gamma$, then the following result holds in the
sense of formal power series \begin{equation}\left\langle
\dprod_{i=1}^n X(x_i,\,\eta)\right\rangle =\ds_{m=0}^\infty
\lambda^m \ds_{{{j_1,\,...,j_n \geq 0} \atop {j_1+...j_n=m}}\atop
{{T_1 \in {\cal A}(j_1)},\,...,\,T_n \in {\cal A}(j_n)}} \Big\langle
{\cal B}_1(T_1,\,x_1,\,\eta)\cdot\cdot\cdot{\cal
B}_n(T_n,\,x_n,\,\eta) \Big\rangle\end{equation}
\end{prop}
\noindent{\bf Proof.}\\ \begin{eqnarray} \label{1.7eqa} \left\langle
\dprod_{i=1}^n X(x_i,\,\eta)\right\rangle &=& \left\langle
\ds_{j=0}^\infty \lambda^j
X_j(x_1,\,\eta)\cdot\cdot\cdot\ds_{j=0}^\infty \lambda^j
X_j(x_n,\,\eta) \right\rangle \,\, \nonumber\\&=& \ds_{m=0}^\infty
\lambda^m \ds_{{{j_1,\,...,j_n \geq 0} \atop
{j_1+...+j_n=m}}}\left\langle X_{j_1}(x_1,\,\eta)\cdot\cdot\cdot
X_{j_n}(x_n,\,\eta)\right\rangle \nonumber\\&=&\ds_{m=0}^\infty
\lambda^m \ds_{{{j_1,\,...,j_n \geq 0} \atop {j_1+...+j_n=m}}\atop
{{T_1 \in {\cal A}(j_1)},\,...,\,T_n \in {\cal A}(j_n)}}
\left\langle {\cal B}_1(T_1,\,x_1,\,\eta)\cdot\cdot\cdot{\cal
B}_n(T_n,\,x_n,\,\eta) \right\rangle
\end{eqnarray}
The last equality is an immediate consequence of Theorem \ref{1thm1}.\\ {\hspace*{15cm \prend}}\\
Let $J\subset \N$ be a finite set. The collection of all partition
of $J$ is denoted by ${\cal D}(J)$. A partition is a decomposition
of $J$ into disjoint, nonempty subsets, i.e
\begin{equation} 
I \in {\cal
D}(J)\Longleftrightarrow \exists \,\,k \in \N,
\,I=\{I_1,...,I_k\},\,I_j\cap I_l=\emptyset\,\, \forall\,\, 1\leq
j<l\leq k\,,\,\cup_{l=1}^k I_l=J.
\end{equation}
\begin{Def}
{\rm 
Let $z_1,...,z_k \in \Gamma$, $I$ a partition of the set
$\{1,...,n\}$, $I \in {\cal D}(\{1,...,n\}),\,\\I={\{I_1,...,I_k\}}$
the truncated moments functions $\langle
\eta(z_1)\cdot\cdot\cdot\eta(z_n)\rangle^{T}$ are recursively
defined by :\begin{equation} \langle \dprod_{i=1}^n \eta(z_i)
\rangle=\ds_{I \in {\cal D}(\{1,...,n\}) \atop I=\{I_1,...,I_k\}}
\dprod_{l=1}^k \langle I_l\rangle^T\end{equation}
where $\langle I_l\rangle^T =\langle \dprod_{j \in
I_l}\eta(z_j)\rangle^T$.
}\end{Def}
\begin{Def}{\label{1.def}} 
{\rm
For ${\cal{G}} \in
P(m,\,n)$ and $x_1,\ldots,x_n\in\Gamma$ we define the number ${\cal
P}[{\cal{G}}](x_1,...,x_n)$ as follows:\\
1- Assign  values  $x_1,...,x_n \in \Gamma$ to the roots
of the trees $T_1,\,...,\,T_n$.\\
\hspace*{2mm}-Assign  values $y_1,...,y_m \in \Gamma$, to the inner
vertices.\\
\hspace*{2mm}-Assign values $z_1,...,z_k \in \Gamma\,,\, $ to the
leaves of type one and assign  values ${\bf z}'_1,...,{\bf z}'_l \in
L_\delta$ to the leaves of type two. \\
2-For every edge in a tree $T_r$, $r=1,\ldots,n$ with two end points ,$e=\{v,\,w\}\,,\,v \leq w$,
assign value $G(e)$ to this edge.\\
3-For each inner empty vertices with $l$ legs connected to leaves with arguments $q_1,\ldots,q_l$ multiply with
\begin{equation}
\langle \eta(q_1) \cdot \cdot \cdot \eta(q_l) \rangle^T.\end{equation} 
4-For the $j$-th leaf of type two multiply with $f({\bf z}'_j)$.\\
5-For each vertex $v \in V({\cal{G}})$ that is an inner vertex of a tree multiply with the
corresponding
multiplicity $M(v)$,\, (see Def.\ref{def11}). \\
6-Integrate with respect to the Lebesgue measure
$dy_1\cdot\cdot\cdot dy_m\,dz_1\cdot\cdot\cdot dz_k\,d{\bf
z}'_1\cdot\cdot\cdot d{\bf z}'_l.$
}\end{Def}
Let $m,\,n \,\in \,\N$, $I \in {\cal D}(\bigcup_{j=1}^n
L_1(T_j)),\,I={\{I_1,...,I_k\}}$ and ${T_1 \in {\cal
A}(j_1)},\,...,\,T_n \in {\cal A}(j_n)$,\, such that $
j_1+...+j_n=m$. We construct a one to one correspondence between
pairs $(T_1,\,T_2,\,...,\,T_n\,;\,I)$ and the generalized Parisi-Wu
graph ${\cal{G}}\,\in \, P(m,\,n)$ in the following way:

 Let the Parisi-Wu graph 
${\cal{G}}$ is
given and we construct $(T_1,...,T_n\,;\,I)$.
By definition ${\cal{G}}$ contains $n$ rooted trees, ${T_1 \in {\cal
A}(j_1)},\,...,\,T_n \in {\cal A}(j_n)$ such that $j_1+...+j_n
=m$.
Suppose that there are $k$ empty vertices in the graph ${\cal{G}}$.
Give an arbitrary number $l=1,...,k$ to each empty
 vertex.
For the $l$-th empty vertex let $I_l$ be the subset in $
\bigcup_{j=1}^n L_1(T_j)$ corresponding to the leaves connected to this
vertex. Then
$I=\{I_1,...,I_k\} \in {\cal D}(\bigcup_{i=1}^n L_1(T_i))$.

Conversely, let $ T_1,...,T_n\,,\,I$ be given such that ${T_1 \in
{\cal A}(j_1)},\,...,\,T_n \in {\cal A}(j_n)$ and $I \in {\cal
D}(\bigcup_{j=1}^n L_1(T_j)),\linebreak \,I={\{I_1,...,I_k\}}$.
 Let $Q=\{1,...,k\}$ be the set of the empty
vertices of ${\cal{G}}$.
 Let $V({\cal{G}})=\sum_{j=1}^n V(T_j)\cup Q$ be the set of all vertices of ${\cal{G}}$
and $E(T_j,\,Q)$ be the set of all edges which connect the $j-th $
rooted tree $T_j,\,(1\leq j\leq n)$ with the empty vertices in $Q.$
Then $E(T_j,\,Q)=\{(v,\,q)\,:\,v \in \,L_1(T_j),\,q \in \,Q\,;\, v
\in
I_q\}$.
 Let $E({\cal{G}})=\sum_{j=1}^n E(T_j)\cup E(T_j,\,Q)$ be the set of all edge of
${\cal{G}}$.

We have established the following result :
\begin{lem}{\label {1.lem}}Let $m,\,n \in \,\N$ and $M(j)=\{T_1 \in {\cal
A}(j_1),\,...,\,T_n \in {\cal A}(j_n)\,;\,\,I \in {\cal
D}(\sum_{j=1}^n L(T_j)),\,I=\{I_1,...,I_k\} \}.$ Then the mapping
$N$ constructed above \begin{equation}N:P(m,\,n)\longrightarrow
M(j)\end{equation}
\end{lem}
is a one to one correspondence.
\begin{theo}\label{theor1}
Let $\{X_t,\, t \geq 0 \}$, be the solution of the stochastic
differential equation (1), and $T_1,...,T_n$, $n$ rooted trees with
roots $x_1,...,x_n.$ Then the moments of $X(x)$ are given by a sum
over all generalized Parisi-Wu graphs ${\cal{G}} \in P(m,\,n)$ of
$m-th$ order that are evaluated according to the ruled fixed in
definition (\ref {1.def}). i.e :\begin{equation}\left\langle
\dprod_{i=1}^n X(x_i,\,\eta)\right\rangle=\ds_{m=0}^\infty\lambda^m
\ds_{{\cal{G}} \in P(m,\,n)}{\cal
P}[{\cal{G}}](x_1,...,x_n)\end{equation}

\end{theo}
\noindent{\bf Proof.} Let $${L_1}({\cal{G}})=\bi_{i=1}^n
L_1(T_{j_i}),\,{L_2}({\cal{G}})=\bi_{i=1}^n
L_2(T_{j_i}),{V}({\cal{G}})=\bi_{i=1}^n
V(T_{j_i})\setminus(\{x_1,...,x_n\}) ,\,{E}({\cal{G}})=\bi_{i=1}^n
E(T_{j_i})$$ The multiplicity of the graph ${\cal G}$ is $M({\cal
G})\,=\dprod_{v\,\in\,V({\cal G})}\,M(v).$\\we have :
\begin{eqnarray} \left\langle \dprod_{i=1}^n X(x_i,\,\eta)\right\rangle
&=&\ds_{m=0}^\infty \lambda^m \left\langle \ds_{{{j_1,\,...,j_n \geq
0} \atop {j_1+...+j_n=m}}\atop {{T_1 \in {\cal A}(j_1)},\,...,\,T_n
\in {\cal A}(j_n)}} {\cal B}_1(T_1,\,x_1,\,\eta)\cdot\cdot\cdot{\cal
B}_n(T_n,\,x_n,\,\eta) \right\rangle\nonumber\\&=&\ds_{m=0}^\infty
\lambda^m \ds_{{{j_1,\,...,j_n \geq 0} \atop {j_1+...+j_n=m}}\atop
{{T_1 \in {\cal A}(j_1)},\,...,\,T_n \in {\cal A}(j_n)}}M({\cal
G})\int \dprod_{e\in E({\cal{G}})}G(e) \nonumber\\&\times&\dprod_{{
l} \in {L_2}({\cal{G}}) }F_f({l}) \left\langle \dprod_{l_1 \in
{L_1}({\cal{G}})}\eta(l_1)\right\rangle \otimes_{v \in
{V}({\cal{G}})\setminus Q({\cal G})}dv \nonumber\\&=&\ds_{m=0}^\infty \lambda^m
\ds_{{{j_1,\,...,j_n \geq 0} \atop {j_1+...+j_n=m}}\atop {{T_1 \in
{\cal A}(j_1)},\,...,\,T_n \in {\cal A}(j_n)}}M({\cal G}) \int
\dprod_{e \in E({\cal{G}})}G(e) \nonumber\\&\times&\dprod_{{l} \in
{L_2}({\cal{G}})}F_f({l})\ds_{I \in {\cal D}({L_1}({\cal{G}})) \atop
I=\{I_1,...,I_k\}} \dprod_{l=1}^k \left\langle \dprod _{l_1 \in
I_{l}}\eta(l_1)\right\rangle^T \otimes_{v \in {V}({\cal{G}})\setminus Q({\cal G}) }dv
\nonumber\\&=&\ds_{m=0}^\infty \lambda^m \ds_{{{j_1,\,...,j_n \geq
0} \atop {j_1+...j_n=m}}\atop {{{T_1 \in {\cal A}(j_1)},\,...,\,T_n
\in {\cal A}(j_n)}\atop{ I \in {\cal D}({L_1}({\cal{G}})) \atop
I=\{I_1,...,I_k\}}}}M({\cal G})\int \dprod_{e \in
E({\cal{G}})}G(e)\nonumber\\&\times&\dprod_{{l} \in
{L}({\cal{G}})}F_f({ l})\dprod_{l=1}^k \left\langle \dprod _{l_1 \in
I_{l}}\eta(l_1)\right\rangle^T \otimes_{v \in V({\cal{G}})\setminus (Q{\cal G}) }dv.
\end{eqnarray} Now we apply Lemma \ref {1.lem} and Definition \ref
{1.def} to conclude.\\{\hspace*{15cm \prend}} 

We denote the
collections of the connected generalized Parisi-Wu graph of $m$-th
order with $n$ roots by
$P_{c}(m,\,n).$\\
We construct a one to one correspondence between $P(m,\,n)$ and
${\cal N}(m,n)=\{ I \in {\cal D}(\{1
...n\}),\,\,I=\{I_1,\,...,\,I_k\}\,;\,\,G_1 \in
P_{c}(m,\,p_1),...,G_k \in
P_{c}(m_k,\,p_k),\,\,m=m_1+...+m_k,\,\,n=p_1+...+p_k \} $ in the following way:

We first consider the case when a generalized Parisi-Wu graph ${\cal{G}} \in P(m,\,n)$ is given.
By definition ${\cal{G}}$ contains $n$ rooted trees $T_1 \in {\cal
A}(j_1),...,T_n \in {\cal A}(j_n)$, with roots $x_1,...,x_n$ and an
arbitrary number
of empty vertices.
For $1\leq q< j \leq n$ we say that the graph connects $q$ and $j$,
in notation $q \sim_{R}j,$ if $x_q$ and $x_j$ are connected in the
generalized Parisi-Wu graph. Obviously $\sim_R$ is
an equivalence relation on $\{1,...,n\}$. Let $I=\{I_1,...,I_k\}$
be the equivalence classes of $\sim_R$ then $I\in \,{\cal
D}(\{1...n\})$. For $1\leq l\leq k$ let ${\cal G}_l$ be the
connected generalized Parisi-Wu graph with $m_l$ inner vertices and
$p_l=\sharp I_l$ roots, then $m=m_1+...+m_k$ and $n=p_1+...+p_k$.
Conversely let $I \in {\cal D}(\{1...n\}),\,
I=\{I_1,\,...,\,I_k\},\, {\cal G}_1 \in P_{c}(m,\,p_1),...,{\cal
G}_k \in P_{c}(m_k,\,p_k),\,\,\,$\linebreak $m=m_1+...+m_k.$

Obviously the connected generalized Parisi-Wu graphs ${\cal G}_1,...,{\cal G}_k$
determine a generalized Parisi-Wu graph ${\cal{G}}$ with
$m=m_1+...+m_k$ inner vertices and $n=p_1+...+p_k$ roots by taking the disjoint union.

This second operation is the inverse of the first one. We
have established the following result:
\begin{lem}\label{baba}
The map \begin{equation}M: P(m,\,n)\longrightarrow {\cal
N}(m,n)\end{equation}is a one to one correspondence.
\end{lem}
The following is known as linked cluster theorem  and gives the
truncated moments of $X(x)$  as a sum over the connected Parisi-Wu
graphs, only.
\begin{theo}
Let $ P_{c}(m,\,n)$ be the set of connected generalized Parisi-Wu
graphs, then the truncated moments $ \langle
X(x_1,\,\eta)\cdot\cdot\cdot X(x_n,\,\eta \rangle^T$ are given by
\begin{equation} \label{ez0}\langle
X(x_1,\,\eta)\cdot\cdot\cdot X(x_n,\,\eta \rangle^T=\ds_{m=0}^\infty
\lambda^m \ds_{{\cal{G}} \in P_{c}(m,\,n)}{\cal
P}[{\cal{G}}](x_1,...,x_n).\end{equation}
\end{theo}
\noindent{\bf Proof.} By Theorem (\ref{theor1}) we have
\begin{equation}\label{ze4}
\langle X(x_1,\,\eta)\cdot\cdot\cdot X(x_n,\,\eta
\rangle=\ds_{m=0}^\infty\lambda^m \ds_{{\cal{G}} \in P(m,\,n)}{\cal
P}[{\cal{G}}](x_1,...,x_n).\end{equation} In the other hand by
definition of the truncated moments \begin{equation}
\label{ze1}\langle X(x_1,\,\eta)\cdot\cdot\cdot X( x_n,\,\eta)
\rangle= \ds_{I \in {\cal D}(\{1,...,n\}) \atop I=\{I_1,...,I_k\}}
\dprod_{l=1}^k \langle \dprod_{j \in I_l} X(x_j,\,\eta) \rangle^T
\end{equation}
To prove the theorem it suffices to prove that in equation
(\ref{ze1}) we can replace the truncated moments by the right hand
side of (\ref{ez0}), i.e.
\begin{equation} \label{ze2}\langle
X(x_1,\,\eta)\cdot\cdot\cdot X(x_n,\,\eta) \rangle= \ds_{I \in {\cal
D}(\{1,...,n\}) \atop I=\{I_1,...,I_k\}} \dprod_{l=1}^k
{\left(\ds_{m_l=0}^ \infty {\lambda^{m_l}} \ds_{{\cal G}_l \in P_{c}
(m_l,\,\sharp I_l)}{\cal P}[{\cal G}_l](x_j\,;\,j \in I_l) \right)}.
\end{equation}
Thus we have to prove that the right hand side of $(\ref{ze2})$ is
equal to the right hand side (r.h.s) of (\ref{ze4}). We have for the
r.h.s of (\ref{ze2})
\begin{eqnarray}
&&\ds_{I \in {\cal D}(\{1,...,n\}) \atop
I=\{I_1,...,I_k\}}\ds_{m_1=0 \atop{\cal G}_1 \in PW_{c}
(m_1,\,p_1)}^\infty\lambda^{m_1}\cdot\cdot\cdot\ds_{m_k=0 \atop
{\cal G}_k \in P_{c} (m_k,\,p_k)}^\infty\lambda^{m_k} \,\,{\cal
P}[{\cal G}_1](x_j\,;\,j \in I_1)\cdot\cdot\cdot{\cal P}[{\cal
G}_k](x_j\,;\,j \in I_k)\nonumber\\&=&\ds_{I \in {\cal
D}(\{1,...,n\}) \atop I=\{I_1,...,I_k\}}\ds_{m_1,...,m_k=0 \atop
{{\cal G}_1 \in P_{c}
(m_1,\,p_1)\atop{\hspace*{0.3cm}.\atop{\hspace*{0.3cm}.\atop{\hspace*{0.3cm}.\atop{\hspace*{0.3cm}{{{\cal
G}_k \in P_{c} (m_k,\,p_k)}}}}}}}}^\infty\lambda^{m_1+...+m_k} {\cal
P}[{\cal G}_1](x_j\,;\,j \in I_1)\cdot\cdot\cdot{\cal P}[{\cal
G}_k](x_j\,;\,j \in I_k)\nonumber\\&=&\ds_{m_1,...,m_k=0}^\infty
\lambda^{m_1+...+m_k}\ds_{I \in {\cal D}(\{1,...,n\}) \atop{
I=\{I_1,...,I_k\} \atop {{\cal G}_1 \in P_{c}
(m_1,\,p_1)\atop{{\hspace*{0.3cm}.\atop{\hspace*{0.3cm}.\atop{\hspace*{0.3cm}.\atop{\hspace*{0.3cm}
 {{\cal G}_k \in P_{c}
(m_k,\,p_k)}}}}}}}}}{\cal P}[{\cal G}_1](x_j\,;\,j \in
I_1)\cdot\cdot\cdot{\cal P}[{\cal G}_k](x_j\,;\,j \in
I_k)\nonumber\\&=&\ds_{m=0}^\infty \lambda^m \,\ds_{I \in {\cal
D}(\{1,...,n\}) \atop{ I=\{I_1,...,I_k\} \atop {{\cal G}_1 \in P_{c}
(m_1,\,p_1)\atop{\atop{{\hspace*{0.3cm}.\atop{\hspace*{0.3cm}.\atop{\hspace*{0.3cm}.\atop{\hspace*{0.3cm}
{{\cal G}_k \in P_{c}
(m_k,\,p_k)\atop{m_1,...,m_k=0,1,...\atop{m_1+...+m_k=m}}}}}}}}}}}}{\cal
P}[{\cal G}_1](x_j\,;\,j \in I_1)\cdot\cdot\cdot{\cal P}[{\cal
G}_k](x_j\,;\,j \in I_k)
\end{eqnarray} Now application of Lemma (\ref{baba}) concludes the argument.
\prend
\section{L\'evy noise }
In this section we
define first a white noise measure with L\'evy
characteristic $\psi$, then we recall a theorem which give the
truncated moments of the noise which permit us to derive a simplification
of the Feynman rules (Definition \ref{1.def}).\\
Let $\nu$ be an infinitely divisible probability distribution. By Levy-Khinchine theorem
(see e.g. \cite{DA}) we know that the Fourier transform (or
characteristic function) of $\nu$, denoted by $C_{\nu}$, satisfies
the following formula
$$ C_{\nu}(t)=\int_{\R} e^{ist}d{\nu}(s)=e^{\psi(t)},\,\,t \in \R$$
where $\psi : \R \longrightarrow \C$ is a continuous function,
called the L\'evy characteristic of $\nu$, which is uniquely
represented as follows
\begin{equation} \label{1.a eqa} \psi(t)=i\bar at-{{\sigma^2 t^2}\over 2}+\int_{\R\setminus
\{0\}}\left(e^{ist}-1-\frac{ist}{1+s^2}\right) \, dM(s),\,\,\,\,\, \forall \,\,t \in\,\, \R.
\end{equation}
where $\bar a \,\in \R$, $\sigma^2\geq 0$ and $M$ is a L\'evy
measure on $\R\setminus\{0\}$, i.e. $\int_{\R\setminus\{0\}}\min(1,s^2)dM(s)<\infty$.
Under the additional assumption that all moments of $M$ exist, i.e. $\int_{\R\setminus\{0\}}s^{2n}dM(s)<\infty$, $n\in\N_0$, one can reparametrize
this representation setting $z=\int_{\R\setminus\{0\}}dM(s)$, $r=M/z$, $a=\bar a-\int_{\R\setminus\{0\}}\frac{s}{1+s^2}dM(s)$ and one obtains\\
\begin{equation} \label{1.aa eqa} \psi(t)=iat-{{\sigma^2 t^2}\over 2}+z\int_{\R\setminus
\{0\}}\left(e^{ist}-1\right) \, dr(s),\,\,\,\,\, \forall \,\,t \in\,\, \R.
\end{equation}
The first term in ({\ref{1.aa eqa}}) is called deterministic, the
second one Gaussian term and the third one Poisson term. \\
We note by ${\cal B}$ the $\sigma$-algebra generated by the cylinder
sets of $S'(\Gamma)$. Then $S'(\Gamma)$ is a measurable space.\\
We define a characteristic functional on $S(\Gamma)$, as a
functional $C:S(\Gamma)\,\longrightarrow\,\C$ such that :\\
\begin{enumerate}\item $C$ is continuous on $S(\Gamma)$;\\
\item $C$ is positive-definite;\\
\item $C(0)=1.$
\end{enumerate}
By the well-known Bochner-Minlos theorem (see \cite{GV}) there
exists a one to one correspondence between characteristic functional
$C$ and probability measures $\mu$ on $(S'(\Gamma),\,{\cal B})$
given by the following relation
\begin{equation}C(f)=\int_{S'(\Gamma)} e^{i\langle f,\,\xi\,\rangle}
\,d\mu(\xi),\, \,\,f \,\,\in \,S(\Gamma).\end{equation} We set
$\Gamma_\pm=\{x\in\Gamma:\, x=(t,{\bf x}),\, \pm t>0\}$.
\begin{theo}\label{1.t1thm}
Let $\psi$ be a L\'evy characteristic given by the
representation (\ref{1.a eqa}) then there exist an unique
probability measure $P_{\psi}$ on $(S'(\Gamma),\,{\cal B})$ such
that 
\begin{equation}\label{1.00eqa}
C_{\eta}(f)=\int_{S'(\Gamma)}e^{i \eta( f,\,\xi)
}dP_{\psi}(\xi)=\exp\left( \int_{\Gamma_+} \psi(f(x))\,dx
\right),\,\,\,\,x=(t,{\bf x}).
\end{equation}
Where \,\,$\eta(f,\,\xi)=\langle\,f,\,\xi \rangle,\,\,\,\,f \in
S(\Gamma).$
\end{theo}
\noindent{\bf Proof.} The right-hand side of (\ref{1.00eqa}) is a
characteristic functional on $S(\Gamma)$ (see eg. \cite{GV}). The
result thus holds by using the Bochner-Minlos Theorem.\\ {\hspace*{15cm
\prend}}
\begin{Def}
{\rm We call $P_{\psi}$ in Theorem \ref{1.t1thm} a generalized white
noise measure with L\'evy characteristic $\psi$ and
$(S'(\Gamma),\,{\cal B},\, P_{\psi})$ the generalized white noise probability
space associated with $\psi$. The associated coordinate process
$$\eta :S(\Gamma)\times (S'(\Gamma),\,{\cal B},\,P_{\psi})\longrightarrow\,\C,~~\eta(f)(\omega)=\omega(f)~\forall f\in S(\Gamma),~ \omega\in S'(\Gamma)$$
is called a L\'evy noise.
}\end{Def}

The following result on L\'evy noise is essential for our perturbative approach

\begin{prop}\label{PROP} There exists a version of $G*\eta$ such that
$G*\eta\in {\cal M}_b$ $P_\psi$ a.s. .
\end{prop}
\noindent {\bf Proof.} The proof requires several steps.

i) By (\ref{1.aa eqa}), the L\'evy noise can be decomposed into
independent parts  $\eta=\eta^d+\eta^g+\eta^p$, where
$\eta^d,\,\,\eta^g$ and $\eta^p$ are respectively the deterministic,
Gaussian and Poisson noise given by their characteristic function in
theorem (\ref{1.t1thm}). It is sufficient to deal with these parts
separately. The statement holds trivially for the deterministic
part.

Let us next focus on the Gaussian part. We get the following two results that together show $G*\eta^g\in{\cal M}_b$ a.s. .

ii) $G\ast \eta^g(t,\,{\bf x})$, for ${\bf x}\in L_\delta$ fixed, has a continuous extension in $t$. In particular, $G*\eta^g$ is a measurable function on $\Gamma$ (a.s.).\\
By the Linked Cluster theorem we have
\begin{equation}
\langle \mid G\ast \eta^g(t,{\bf x})-G\ast \eta^g(s,{\bf x})\mid^4\rangle=3(\langle
\mid G\ast \eta^g(t,{\bf x})-G\ast \eta^g(s,{\bf x})\mid^2\rangle)^2.
\end{equation}
If we assume that $t>s$, then it is easy to see, using $\langle\eta^g(t,{\bf x})\eta^g(s,{\bf y})\rangle=\sigma^2\theta(t)\theta(s)\delta(t-s)\delta({\bf x}-{\bf y})$, that
\begin{eqnarray}\label{array} &&\left\langle \mid G\ast \eta^g(t,{\bf x})-G\ast
\eta^g(s,{\bf x})\mid^2\right\rangle\nonumber\\
&&~~~~~~~~~~~~=\left\langle (G\ast
\eta^g(t,{\bf x}))^2+(G\ast \eta^g(s,{\bf x}))^2-2(G\ast \eta^g(t,{\bf x}))(G\ast
\eta^g(s,{\bf x}))\right\rangle\nonumber\\
&&~~~~~~~~~~~~=\Bigg\langle \Big(
\int_{0}^\infty\int
\theta(t-t')\tilde{G}_{t-t'}({\bf x}-{\bf x}')\eta^g(t',{\bf x'})d{\bf x}'dt'\Big)^2\nonumber\\
&&~~~~~~~~~~~~~~~~~~~~~~~~+\Big(
\int_{0}^\infty\int
\theta(s-t')\tilde{G}_{s-t'}({\bf x}-{\bf x}')\eta^g(t',{\bf x}')d{\bf x}'dt'\Big)^2\nonumber\\
&&~~~~~~~~~~~~~~~~~~~~~~~~-2\,\int_0^\infty
\int_0^\infty\int\int\theta(t-t')
\theta(s-s')\tilde{G}_{t-t'}({\bf x}-{\bf x}')\tilde{G}_{s-s'}({\bf x}-{\bf x}'')\nonumber\\
&&~~~~~~~~~~~~~~~~~~~~~~~~~~~~~~~~~~~~~~~~~~~~~~~~~~~~~~~~~~~~~~~~~~~~~~\times\eta^g(t',{\bf x}')\eta^g(s',{\bf x}'')d{\bf x}d{\bf x}''dt'ds'\Bigg\rangle
\nonumber\\
&&~~~~~~~~~~~=\sigma^2\, \int_0^\infty\int
\Big(\theta(t-t')\tilde{G}_{t-t'}^2({\bf x}-{\bf x}')+\theta(s-t')\tilde{G}_{s-t'}^2({\bf x}-{\bf x}')\nonumber\\
&&~~~~~~~~~~~~~~~~~~~~~~~~-2\theta(s-t')\tilde{G}_{t-t'}({\bf x}-{\bf x}')\tilde{G}_{s-t'}({\bf x}-{\bf x}')
\Big)d{\bf x}'dt'\nonumber\\
&&~~~~~~~~~~~=\sigma^2\,\int_0^\infty\int \theta(s-t')
\Big(\tilde{G}_{t-t'}({\bf x}')-\tilde{G}_{s-t'}({\bf x}')\Big)^2d{\bf x}'dt'+\int_s^t\int\tilde{G}_{t-t'}^2({\bf x}')d{\bf x}  'dt'
\end{eqnarray}
Moreover by the expression of $\tilde{G}_t$ given in
$(\ref{1.12eqa})$ we get :
\begin{equation}
\left|(1+|{\bf x}|^2)^{d}\left(\tilde{G}_{t}({\bf x})-\tilde{G}_{s}({\bf x})\right)\right|^2 \leq {{c^2} (t-s)^2\, \, \over
{(2\pi)^{2d}}}~\mbox{and}~\left|(1+|{\bf x}|^2)^dG_{t-t'}^2({\bf x})\right|\leq c\,\frac{\delta^{-2\,d}}{(2\pi)^{2d}}.
\end{equation}
Here $c=\sup_{t>0}\int_{\Pi_{\delta}^d}(1-\Delta)^de^{-t\mu_{\delta,\,m}^2({ \bf p})}d{\bf
p}<\infty$. Inserting these two estimates in (\ref{array}),  we get for $K>0$ sufficiently large \begin{equation} \langle \mid G\ast \eta^g(t,{\bf x})-G\ast
\eta^g(s,{\bf x})\mid^4 \rangle \leq K \mid t-s\mid^{2}.\end{equation}
Application of Kolmogorov's extension theorem now proves the assertion.

ii) $G\ast \eta^g(t,\,{\bf x})$ is polynomially bounded (a.s.).\\
We prove the stronger statement that the expectation of $\int{
(G\ast \eta^g( x))^k\over{ (1+\mid x \mid^2)^N}}\, dx$ is finite. We
have for $k\in\N$ even
\begin{eqnarray}\left \langle\int{
(G\ast \eta^g( x))^k\over{ (1+|x |^2)^N}}\, dx\right\rangle
&=&A\int {\Big(\left\langle \mid G\ast \eta^g({\bf x})\mid^2
\right\rangle \Big)^{k\over 2} \over {(1+\mid
x\mid^2)^N}}\, dx\nonumber\\
\nonumber\\&\leq& A\sigma^2\left(\int_{\Gamma}G^2(x)\,dx\right)^{\frac{k}{2}}\int {dx \over{ (1+\mid x\mid^2)^N}} <\infty.
\end{eqnarray} The last inequality hold for an arbitrary
$N\,>d/2.$ Here $A$ is the number of pairings of $k$ objects.

It is left over to deal with the poisson part. First we recall the path properties of $\eta^p$ following \cite{AGY}.

iii) $\eta^p$ is a signed measure with locally discrete support (a.s.). \\
Let $\Lambda_n\,\subseteq\,\Gamma_+$ be a monotone sequence
of compact sets s.t $ \Lambda_n\,\uparrow\,\Gamma_+$ as
$n\,\rightarrow\,\infty$ and $\Lambda_0=\emptyset.$ For
$n\,\in\,\N$ let $D_n=\Lambda_n\,\setminus\,\Lambda_{n-1}$ and we
denote the (Lebesgue) volume of $D_n$ by $\mid\,D_n\mid$. Let
$\tilde\eta_n^p =\sum_{j=1}^{N^z_n}S^n_j\delta_{X^n_j}$ be a random field, where $\delta_{x}$ is the Dirac measure of mass one in $x$ and
$N^z_n$ is a
Poisson random variable  with intensity
$z\mid\,D_n\mid$ i.e
$$P(N^z_n=l)=e^{-{z\mid\,D_n\mid}}{(z\mid\,D_n\mid)^l \over
{l!}}\,\,;\,\,l \,\in\,\N_0.$$
For $n\in\N$, $\{X_j^n\}_{j\in\N}$ is a family of i.i.d. $\Gamma$-valued random variables distributed uniformly on $D_n$. $\{S_j^n\}_{j,n\in\N}$ is a family of real valued random variables
with law given by $r$. All these random variables are independent of each other.
The characteristic functional of
${\tilde\eta_n^p}$ is given by
\begin{eqnarray}\left\langle e^{i\eta_n^p(f)}\right\rangle&=&\left\langle e^{i\sum_{j=1}^{N^2_n}S^n_j f({X^n_j})}\right\rangle
\nonumber\\&=&e^{-z\mid\,D_n\mid}\sum_{l=0}^\infty\,{(z D_n)^l \over
{l!}}\left(\int_{D_n}\int_{\R\setminus\{0\}}
e^{isf(x)}\, dr(s)\frac{dx}{|D_n|}\right)^l\nonumber\\&=&\exp\left\{z\int_{D_n}\int_{\Gamma
\backslash\{0\}}(e^{is\,f(x)}-1)\,dr(s)dx\right\}={\cal C}_{\eta^p}(f), ~~~~~\forall f\in S(\Gamma), {\rm supp}f\subseteq D_n.\nonumber\\
\end{eqnarray}
Hence $\eta^p$, when restricted to $D_n$, coincides in law with
$\eta^p_n$. By the uniqueness statement of the Bochner-Minlos
Theorem and the fact that by construction $\tilde\eta_n^p$ takes
values in the locally discrete signed measures (a.s.), it follows
(cf. \cite{AGY} for the details) that also $\eta^p$ with probability
one is a locally discrete signed measure.

iv) For $f\in L^1(\Gamma,dx)$ bounded and non-negative, $\eta^p$ is $f$-finite (a.s.), i.e. $\int_\Gamma f\, d|\eta^p|<\infty$. Here $|\eta^p|=\eta^p_++\eta^p_-$ is the modulus of the signed measure $\eta$, cf. \cite{Hal}.\\
This assertion can be seen from the following calculation
\begin{eqnarray}\left\langle\int_\Gamma f\,d\mid \eta^p \mid\right\rangle&=&\ds_{n
=0}^\infty\,\left\langle\int_{D_n}\,f\,\,d\mid \tilde\eta_n^p
\mid\right\rangle \nonumber\\&=&\ds_{n =0}^\infty
\left\langle\,\ds_{j=1}^{N_n^z}\,\mid\,S_j^n\mid
f(X_j^n)\right\rangle\nonumber\\&=&\ds_{n=0}^\infty e^{-z\mid
D_n\mid}\,\ds_{l=1}^\infty{(z \mid D_n\mid)^l \over {(l-1)!}}
\int_{D_n}f(x){dx\,\over \mid D_n\mid}\int_{\R\setminus\{0\}}\mid
s\mid\,dr(s)\nonumber\\&=&z\ds_{n=0}^\infty
e^{-z\mid\,D_n\mid}\,e^{z\mid\,D_n\mid}\,\int_{D_n}
f(x)dx\int_{\R\setminus\{0\}} \mid s\mid
dr(s)\nonumber\\&=&z\,\int_{\R\setminus\{0\}} \mid s\mid
dr(s)\,\int_{\Gamma}\,f(x)\,{dx}<\infty\,.
\end{eqnarray}

v) $G\ast{\eta^p}\,\,\in\,\,L^1(\Gamma,\,g_{\epsilon}\,\,dx)$ a.s., where $g_\epsilon=(1+|x|^2)^{-(d/2+\epsilon)}$. In particular, $G*\eta^p$ is measurable.\\
 Let $D_n^l=\Lambda_n\,\setminus\,\Lambda_l$ for $n\,>\,l$, we denote the
restriction of the noise $\eta^p$ to an open set
$A\,\subseteq\,\Gamma$ by $\eta^p_{\mid_A}$.
Clearly,
$G\ast\,\eta^p_{\mid_{\Lambda_n}}\,\in\,L^1(\Gamma,\,g_{\epsilon}\,dx)$ (a.s.)
since $G$ is in $L^1(\Gamma,\,dx)$ and
supp$\eta^p_{\mid_{\Lambda_n}}$ is finite (a.s.). The following estimate
shows that $G\ast\,\eta^p_{\mid_{\Lambda_n}}$ forms a Cauchy
sequence in $L^1(\Gamma,\,g_{\epsilon}\,dx)$. With $\parallel
\,.\,\parallel_{\epsilon,\,1}$ the $L^1-$norm on that space, we get
\begin{eqnarray}
\sup_{n>\,l}\parallel
\,G\ast\,\eta^p_{\mid_{\Lambda_n}}-\,G\ast\,\eta^p_{\mid_{\Lambda_l}}\parallel_{\epsilon,\,1}&=&\sup_{n>\,l}\parallel
\,G\ast\,\eta^p_{\mid_{D_n^l}}\parallel_{\epsilon,\,1}\nonumber\\&\leq&\sup_{n>\,l}\int_{\Gamma}\,\mid
G\mid\ast\,g_{\epsilon}d\mid
\,\eta^p_{\mid_{D_n^l}}\,\mid\nonumber\\&=&\int_{\Gamma \setminus
\Lambda_l}\,\mid G\mid\ast\,g_{\epsilon}d\mid
\,\eta^p\,\mid\,\longrightarrow\,0\,\,as\,\,l\longrightarrow\,\infty
\end{eqnarray}
since $\eta^p$ is $|G\mid\ast\,g_{\epsilon}$-finite (see iv) and note that $L^1(\Gamma,dx)$ is closed under convolutions). Also,
\begin{equation}
\lim_{n\,\longrightarrow\,\infty}\langle\,G\ast\,\eta^p_{{\Lambda_n}},\,f\rangle=\lim_{n\,\longrightarrow\,\infty}
\langle\,\eta^p_{{\Lambda_n}},\,G\ast\,f\rangle=\langle\,\eta^p,\,G\ast\,f\rangle
=\langle\,G\ast\eta^p,\,f\rangle\,\,\,\forall\,f\,\in\,{\cal S},
\end{equation}
and by the fact that convergence in $L^1(\Gamma,\,g_{\epsilon}\,dx)$
implies convergence in ${\cal S'},$ we get that $G\ast\,{\eta^p}$
coincide with the limit of $G\ast\,{\eta^p_{\Lambda_n}}$ in the
Banach space $L^1(\Gamma,\,g_{\epsilon}\,dx).$\\

vi) $G\ast \eta^p(t,\,{\bf x})$ is polynomially bounded.\\
\begin{eqnarray}\int {\left \langle (G\ast \eta^p(
x))^n\right\rangle\over{ (1+\mid x \mid^2)^N}} dx&=&\int\ds_{I \in
{\cal D}(\{1,...,n\}) \atop I=\{I_1,...,I_k\}} \dprod_{l=1}^k
{{\left\langle (G\ast \eta^p(
x))^{\sharp\,I_{l}}\right\rangle^T\over{ (1+\mid x
\mid^2)^N}}}\,dx\nonumber\\&=&\ds_{I \in {\cal D}(\{1,...,n\}) \atop
I=\{I_1,...,I_k\}} \int\dprod_{l=1}^k{\left\langle (G\ast \eta^p(
x))^{\sharp\,I_{l}}\right\rangle^T\over{ (1+\mid x
\mid^2)^N}}\,dx\nonumber\\&\leq&C\int\,{dx\, \over{ (1+\mid
x\mid^2)^N}}<\infty.
\end{eqnarray}
The last inequality hold for an arbitrary $N\,>d/2.$ Here
$$C=b_n\max_{n_l:\sum_{l=1}^k n_l=n}\sup_{(t,{\bf x})\in\Gamma}\prod_{l=1}^n\left|\langle(G*\eta^p(t,{\bf x}))^{n_l}\rangle^T\right|$$
with $b_n$ the $n$-th Bell number, i.e. the number of partitions of $n$ objects, and the supremum of the truncated expectation values is finite by virtue of Theorem \ref{1.2t thm} below and Lemma \ref{l.1} \prend

\begin{theo}\label{1.2t thm}
The truncated moment functions of the L\'evy noise $\eta$ are given by the following formula
\begin{equation}\label{0.0eqa}
\langle \eta(x_1)\cdot\cdot\cdot\,\eta(x_n) \rangle^T
=c_n\,\int_{\Gamma_+} \delta(x-x_1)\cdot\cdot\cdot\,\delta(x-x_n)
dx\,.
\end{equation}
where
\begin{eqnarray} c_n&=&(-i)^n {{{d^n\psi(t)} \over dt^n}\mid t=0
}\nonumber\\&=&\delta_{n,\,1} a+
\delta_{n,\,2}\sigma^2+z\int_{\R\setminus\{0\}}s^n
dr(s)\end{eqnarray} $\delta_{n,\,n'}$ being the Kronecker symbol.
\end{theo}
\noindent{\bf Proof.}  After application of the linked cluster
theorem to (\ref{1.00eqa}), the formula follows by a straight
forward calculation, cf. \cite{AGW} for details. \prend

\begin{center}
\begin{figure}
\psset{xunit=1cm,yunit=1cm,runit=1cm,shortput=tab}
\begin{pspicture}(-1,-1)(1,1) \cnode(3.3,0){5pt}{A_1}
 \cnode*(2,0){2.5pt}{C_1} %\ncline{A_1}{C_1}
\dotnode[dotstyle=x,dotscale=2](1.2,0){D}\ncline{C_1}{D}
%\dotnode[dotstyle=otimes,dotscale=3](5,5){E}
\dotnode[dotstyle=otimes,dotscale=2](2.5,0.5){F_1}
\dotnode[dotstyle=pentagon,dotscale=2](2.8,0){F1}
\ncline{A_1}{F1}\ncline{C_1}{F1}
\dotnode[dotstyle=pentagon,dotscale=2](3.8,0){F2}\ncline{B_2}{F2}\ncline{A_1}{F2}
\dotnode[dotstyle=x,dotscale=2](5.2,0){D_1}\ncline{D_1}{F2} %\ncline{A_1}{D_1}
\dotnode[dotstyle=otimes,dotscale=2](2.5,-0.5){F_2}
\dotnode[dotstyle=otimes,dotscale=2](4.,0.5){B}
\dotnode[dotstyle=otimes,dotscale=2](4.,-0.5){B_1}
\cnode*(4.5,0){2.5pt}{B_2} \ncline{B}{B_2}\ncline{B_2}{B_1}
\ncline{C_1}{F_1} \ncline{C_1}{F_2}
\end{pspicture}
\begin{pspicture}(-1,-1)(1,1)\hspace*{4cm}$\equiv$\hspace*{0.1cm} \cnode(3,0){5pt}{A_1}
 \cnode*(2.,0){2pt}{C_1} \ncline{A_1}{C_1}
\dotnode[dotstyle=x,dotscale=2](1.2,0){D}\ncline{C_1}{D}
%\dotnode[dotstyle=pentagon,dotscale=3](5,5){E}
\dotnode[dotstyle=otimes,dotscale=2](2.5,0.5){F_1}
%\dotnode[dotstyle=square*,dotscale=2](3,0){F1}
%\ncline{A_1}{F1}\ncline{C_1}{F1}
%\dotnode[dotstyle=square*,dotscale=2](5,0){F2}\ncline{B_2}{F2}\ncline{A_1}{F2}
\dotnode[dotstyle=x,dotscale=2](5,0){D_1}\ncline{A_1}{D_1}
\dotnode[dotstyle=otimes,dotscale=2](2.5,-0.5){F_2}
\dotnode[dotstyle=otimes,dotscale=2](3.5,0.5){B}
\dotnode[dotstyle=otimes,dotscale=2](3.5,-0.5){B_1}
\cnode*(4,0){2.5pt}{B_2} \ncline{B}{B_2}\ncline{B_2}{B_1}
\ncline{C_1}{F_1} \ncline{C_1}{F_2}
\end{pspicture}\caption{\label{bg.fig}
Construction of a generalized Feynman graph where every leave of
type one together with the two edges connected to it is replaced by
one edge}
\end{figure}
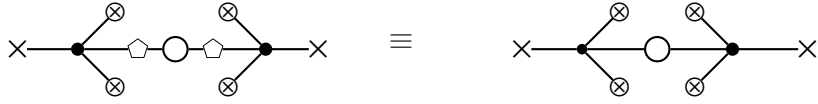
\end{center}

The particular form of (\ref{0.0eqa}) makes it trivial to carry out
the integrals over the leaves of type one of the Parisi-Wu graph.
This leads to the following simplification of Feynman rules:

\begin{theo}\label{mama}
Suppose that the noise $\eta$ is of L\'evy type. Note that in a
Parisi-Wu graph, each leave of type one has exactly two legs. For a
given Parisi-Wu graph ${\cal G}\in P(m,n)$, one can thus get an
equivalent graph ${\cal G}'$ where every leave of type one together
with the two edges connected to it is replaced by one edge (cf.
Figure (\ref{bg.fig})).

 Let $P_c'(m,n)$ be the collection of all
connected graphs obtained in this way from $P_c(m,n)$. The Feynman
rules ${\cal V}[{\cal G}']$ for the graph ${\cal G}'$ then simplify
in the following way: \\
1- Associate the values $x_j$, $j=1,\ldots,n$, to
the roots, an integration variable $y$ to each inner vertex and an
integration variable $z$ to each leave of type 2. \\
2- For each inner vertex $v$ of a tree $T_j$, multiply with a multiplicity factor $M(v)$ and for each
inner empty vertex with $l$ legs multiply with a factor $c_l$.\\
3- For
each leave of type two multiply with a factor $F_f(z)$.\\
4- For each
edge multiply with a propagator $G(e)$. \\
5- Then integrate over all inner vertices
and all leaves of type two. 

In this way, one obtains in the sense of
formal power series
\begin{equation}
\langle X(x_1)\ldots
X(x_n)\rangle^T=\sum_{m=0}^\infty(-\lambda)^m\sum_{{\cal G}'\in
P'_c(m,n)}{\cal V}[{\cal G}'](x_1,\ldots,x_n)
\end{equation}
\end{theo}

\section{Equilibrium correlation function} 

In this short section we determine the limit of the correlation functions when all time arguments go to infinity simultaneously.

\begin{Def}
{\rm Let $P_{1,c}'(m,n)$ be the collection of graphs in $P_c'(m,n)$ that
do not posses a leave of type 2. For ${\cal{G}}\,\in P_{1,c}'(n,m)$
the value ${\cal
V_{\infty}}[{\cal{G}}]({\bf x_1},...,\,{\bf x_n})$ is obtained the following Feynman rules:\\
1- Assign values  $(0,{\bf x_1}),...,(0,{\bf x_n}) \in L_{\delta}$ to the
roots
of the trees $T_1,\,...,\,T_n$ .\\
\hspace*{2mm}-Assign  values  $y_1,...,y_m \in \Gamma$, to the inner
vertices of the tree.\\
\hspace*{2mm}-Assign values $z_1,..., z_k \in
\Gamma$ to the inner empty vertices. \\
2- For every vertex,
$v$ in a tree, multiply with the multiplicity coefficient $M(v)$ and for each inner empty vertex witl $l$ legs by $c_l$.
\\
3-For every edge with two end points $e=\{v,\,w\},\,(v\leq w),$
assign a value $G(e)$ . \\ 
%5-For the $l-$inner full vertices multiply with the multiplicity $M$ of the graphs(See.Remark {\ref{1.rem}}) \\
4-Integrate over $\Gamma_{-}$ with respect to the Lebesgue measure
$dy_1\cdot\cdot\cdot dy_m d z_1 \cdot\cdot\cdot dz_k.$
}\end{Def}
\begin{lem}\label{lemma.1}
Let ${\cal{G}} \in P_c'(m,\,n)$  and let $I({\cal
G})=\prod_{e\,\in\,E({\cal G})}G(e)$, then
$\forall\,N\,\in\,\,\N,\,\,\exists\,\,K=K(N,m)$ such that
\begin{equation}\mid I({\cal G})\mid\,\,\leq\,\,{K\over
{(1+\dmax_{v,\,w\,\in\,V({\cal G})}\mid
v-w\mid^2)^N}}\,\,.\end{equation}
\end{lem}
\noindent{\bf Proof.}\,\,Let $v',v''\in V({\cal G})$ such that $\mid v'-v''\mid^2
\,\,=\dmax_{v,\,w\,\in\,V({\cal G})\setminus \{v_0\} \atop v \neq
w}\mid v-w\mid^2$. As ${\cal G}$ is connected, there exist a
walk $W$ from $v'$ to $v''$ and let $q\leq m+1$ be the number of steps from
$v'$ to $v''$.  Let $e_1,\,e_2,\,...,\,e_q$ be the edges of the walk
$W$. Then
\begin{equation}
 I({\cal
G})=\dprod_{e\,\in\,W}G(e)\,\times\,\dprod_{e\in E({\cal G})\setminus W}G(e)
\end{equation}

In $W$, there must be at least one of the $q$ steps that is
$\geq\,{\mid v'-v''\mid\,\over q}\,\,$. By the use of Lemma
\ref{l.1} we find \begin{equation}\mid I({\cal
G})\mid\,\,\leq\,\,{K(N,m)\over {(1+\mid
v'-v''\mid^2)^N}}\,\,.\end{equation}
 \\
{\hspace*{15cm \prend}}

\begin{theo} \label{6.1theo}Let ${\cal{G}}\,\in \,{{W}}_c(m,\,n),$ then the perturbation
series for the truncated moments converges in the sense of formal
power series, when $t$ goes to infinity to :\begin{equation}\langle
X_{\infty}({\bf x_1})\cdot\cdot\cdot X_{\infty}({\bf x_n})
\rangle^T=\lim_{t \,\longrightarrow\, \infty}\langle
X(x_1)\cdot\cdot\cdot X(x_n) \rangle^T\end{equation} where
\begin{equation}\langle X_{\infty}({\bf x_1})\cdot\cdot\cdot
X_{\infty}({\bf x_n}) \rangle^T=\ds_{m=0}^\infty\lambda^m
\ds_{{\cal{G}} \in {P}_{1,c}'(m,\,n)}{\cal
V}_{\infty}[{\cal{G}}]({\bf x_1},...,{\bf x_n}).\end{equation}
\end{theo}
\noindent{\bf Proof.} Let  ${P}_{2,c}'(m,\,n)={P}_c'(m,\,n)\setminus
{P}_{1,c}'(m,\,n))$. We
have
\begin{equation}
\langle X({ x_1})\cdot\cdot\cdot X({ x_n})
\rangle^T=\ds_{m=0}^\infty\lambda^m \ds_{{\cal{G}} \in
{P}_{1,c}'(m,\,n)}{\cal V}[{\cal{G}}]({ x_1},...,{
x_n})+\ds_{m=0}^\infty\lambda^m \ds_{{\cal{G}} \in
{{P}}_{2,c}'(m,\,n)}{\cal V}[{\cal{G}}]({x_1},...,{ x_n})\hspace*{1cm}
\end{equation}
As we take the limit $t\longrightarrow \infty$ in the sense of
formal power series, it suffices to show
\begin{equation}\lim_{t\longrightarrow \infty}{\cal
V}[{\cal{G}}]((t,{\bf x_1}),...,(t,{\bf x_n}))={\cal V}_\infty[{\cal
G}]({\bf x_1},\ldots{\bf x_n}),\, \forall {\cal{G}} \,\in
{P}_{1,c}'(m,\,n)\end{equation} and $\lim_{t\longrightarrow
\infty}{\cal V}[{\cal{G}}]((t,{\bf x_1}),...,(t,{\bf x_n}))=0,\,
\forall {\cal{G}} \,\in {P}_{2,c}'(m,\,n) $
\\Let first ${\cal{G}} \,\in
{{P}}'_{1,c}(m,\,n)$, we have
\begin{equation} \label{1.01eqa} {\cal V}[{\cal{G}}]((t,{\bf
x_1}),...,(t,{\bf x_n}))=M({{\cal{G}}})\int_{\Gamma_{(0,t)}^m}\dprod_{e
\in E({\cal{G}})} G(e)\dso_{v \in
V({\cal{G}})\setminus \{(t,{\bf x_1}),\,...,\,(t,{\bf x_n})\}} dv
\end{equation}
where $\Gamma_{(a,b)}=\{x=(t,{\bf x})\in \Gamma, a<t\leq b\}$. Now
we  transform  $s\to s-t$ for $v=(s,{\bf v})\in V({\cal G})$ and let
$t$ go to infinity.\\ Let
$\Lambda_e=\{e=\{v,\,w\},\,\,{v=(t_v,\,{\bf{v}})},\,\,\,w=(t_w,\,{\bf{w}}),\,\,
(v \leq w)\},$ we get
\begin{eqnarray}\label{eq1}
\lim_{t\,\longrightarrow\,\infty}{\cal
V}[{\cal{G}}]((t,{\bf x_1}),...,(t,{\bf x_n}))&=&M({{\cal{G}}})\lim_{t\to\infty}\left.\int_{\Gamma^m_{(-t,0)}}\dprod_{\Lambda_e}
G(e) \dso_{v \in {V}({\cal{G}})\setminus(\{ {x_1},\,...,\,{x_n}\})}
d{v}\right|_{x_l=(0,{\bf x_l})}\nonumber\\&=&{\cal V}_{\infty}[{\cal{G}}]({\bf x_1},...,{\bf
x_n})
\end{eqnarray}
where the convergence in the last step is due to Lemma \ref{lemma.1}.

Let now ${\cal{G}} \in {{P}}_{2,c}'(m,\,n)$. Then, again by Lemma
\ref{lemma.1}, $I({\cal G})$ is rapidly decreasing in the difference
of the time argument $t$ of the external vertices and the time
argument $0$ of the leave of type two. Consequently, the integral
over $I({\cal G})$ gives no contribution in the limit $t\to\infty$.
\prend

\section{Determination of the law of the noise}
The main subject of this section is to find statistical information
on the law of the noise $\eta$ driving the SPDE (\ref{1.1eqa}) from
empirical data. One approach to achieve this is the least square
method, i.e. to solve the minimization problem for the two point
function of the stationary distribution
\begin{equation}
\label{7.1eqa}\left\{
\begin{array}{ll}
Q(c_1,\,c_2,...)= \int \mid
F_{th}({\bf x},\,c_1,c_2,...)-F_{em}({\bf x})\mid^2  d{\bf x},\\
{\partial Q(c_1,\,c_2,...) \over \partial c_j}=0,\,\,j\,\in\,\N,
\end{array}
\right.
\end{equation}
where $F_{em}$ is the empirical\footnote{In material sciences, the
empirical correlation function $F_{em}(t, {\bf x},\,c_1,c_2,...)$
can be measured by X-ray spectroscopy or sampling methods, cf. e.g. \cite{BS,HHZ}.} correlation
function.
Note that by Lemma (\ref{lemma.1}) the function $F_{th}({\bf
x})=\langle X_{\infty}({\bf 0})X_\infty({\bf x})\rangle^T$ to any
order of perturbation theory is rapidly decreasing in ${\bf x}$ and
thus there should be no problem with the convergence of the integral
in (\ref{7.1eqa}) -- if the modeling is not completely wrong,
$F_{em}$ should also be of fast decay.

In the remainder of the section, we give the solution to
(\ref{7.1eqa}) to the first order in perturbation theory for $p=3$.
To simplify the calculation we note first
\begin{lem}Let the measure $r$ be symmetric, i.e $r(A)=r(-A)\,\,\forall\,\,\, A \,\,\in {\cal
B}.$ Then one can omit all such generalized Parisi-Wu graphs from
the perturbation series that have an empty vertex with an odd number
of legs.
\end{lem}
We assume this symmetry of
$r$ in the following.
\begin{figure}
\centerline{\psset{xunit=1cm,yunit=1cm,runit=1cm,shortput=tab}}
\begin{center}\begin{pspicture}(-1.5,-1.5)(1.5,1.5)$F_{th}( {\bf
x},\,c_2,c_4)= \hspace*{-4mm}$ \cnode(1.5,0){5pt}{A}
\ncarc[arcangle=45]{B}{C}
\dotnode[dotstyle=x,dotscale=2](0.5,0){D}\ncline{A}{D}
\dotnode[dotstyle=x,dotscale=2](2.5,0){D_1}\ncline{A}{D_1}
\end{pspicture}  \begin{pspicture}(-1.5,-1.5)(1.5,1.5)
$+
 $\cnode(2.5,0){5pt}{A}
 \cnode*(1.5,0){3pt}{C} \ncline{A}{C}
\ncarc[arcangle=45]{C}{A} \ncarc[arcangle=-45]{C}{A}
\dotnode[dotstyle=x,dotscale=2](0.5,0){D}\ncline{C}{D}
\dotnode[dotstyle=x,dotscale=2](3.5,0){D_1}\ncline{A}{D_1}
\end{pspicture}\begin{pspicture}(-1.5,-1.5)(1.5,1.5)
$\hspace*{7mm}+\,\,\,o(\lambda^2)$
\end{pspicture}\caption{\label{b.fi} Expression of the correlation function in the
first order expansion.}\end{center}
\end{figure}
On inspection of the fist order solution, we obtain
by the use of Theorem (\ref{mama}): \\
\begin{equation} \label{Dah1}F_{th}({\bf x},\,c_2,c_4)=c_2\,{P}_1({\bf x})+\lambda\,c_4\,P_2({\bf x})+o(\lambda^2)\end{equation}
where $P_{j}({\bf x})={\cal V}_\infty[{\cal G}_j]({\bf 0},{\bf
x})/c_{2j}$, $j=1,2$ and ${\cal G}_j$ the first/second graph in the
above first order expansion (cf. Figure \ref{b.fi}). Note that $P_j$
does not depend on $c_2,c_4$ anymore.\\ Apparently, $Q$ only depends
on $c_2$ and $c_4$. We have to solve the equations $
{\partial\,Q(c_2,\,c_4)\,\over \,{\partial\,c_i}}=0\,,i=2,\,4$
Moreover,
\begin{equation} Q(c_2,\,c_4)=\alpha\,c^2_2+\beta\,c^2_4+2\gamma\,c_2\,c_4-2\,c_2\,a\,-2\,c_4\,b\,+c\end{equation}
where $\alpha=\int P_1^2d{\bf x}$, $\,\beta=\lambda^2\int P_2^2d{\bf
x}$, $\gamma=\lambda\int P_1P_2d{\bf x}$, $\,a=\int P_1F_{em}d{\bf
x}$, $\,b=\lambda\int P_2F_{em}d{\bf x}$ and $c=\int F_{em}^2d{\bf
x}$. We now write down the equations $ {\partial\,Q\,\over
\,{\partial\,c_i}}=0\,,i=2,\,\,4.$ and solve for $c_2$ and $c_4$.\\
One then  obtains the following first order approximation of $c_2$ and $c_4$:
\begin{equation}
c_2={{a\,\beta-\gamma\,b}\over {\,\alpha\,\beta-\gamma^2}}
\end{equation}
and
\begin{equation}
c_4={{\alpha\,b-\gamma\,a}\over{\alpha\,\beta-\gamma^2}}
\end{equation}
$c_2$ gives a measure for strength of the fluctuations of $\eta$.
Let the kurtosis $K$ given by the following expression:
\begin{equation}
K={c_4 \over c^2_2}\end{equation} If $K=0\, $, then there is no
jump, the stochastic dynamic is purely diffusive. If $0< |K|
<<1\,,$ there are
some jumps but the stochastic dynamic is predominantly diffusive.
If $|K|>>1\,,$ then the stochastic dynamic is predominantly ruled by
jumps.

\section{Appendix A}
For the convenience of the reader we collect some graph-theoretic notions
which have been used in this work.

 Let $V$ be a 
finit set and
$$E=\{e=\{v,\,w\}\,;\,\,v,\,w \in V\}.$$
The elements of $V $ are called vertices, which are of different
types in our case we have the roots $ \times$\,, the inner vertices
$\bullet$ ,\,the empty vertices $\circ$ and the leaves of type one
$\otimes$. Such vertices are labeled by arguments in $\Gamma$.
We have another type of vertex, called the leave of type two,
 ~\begin{pspicture}(0,0)(0,0)\psdot[dotstyle=pentagon,dotsize=7pt](0,.1)\end{pspicture}\hspace{.4cm}, which are in $ L_\delta$.\\ The elements of $E$ are called
edges , i.e lines connecting exactly two vertices, and we say that
an edge $e=\{v,w\}$ joins $v$ and $w$.
Thus we can define a graph
$$ {\cal G}={\cal G}(V)$$ with the vertex set $V$ as a family of 
pairings
$$E({\cal G})\subseteq\{e=\{a,\,b\}, \, a,\,b \in V\}.$$
For a given graph ${\cal G}$ we note $E({\cal G})$ (resp. $V({\cal G})=V$) the set of edges (resp. vertices ) of ${\cal G}$.\\
A graph $H$ is called a subgraph of the graph ${\cal G}$ and we
write
$H \subset {\cal G}$ when the vertex set $V(H)$ of $H$ is contained in the vertex set $V({\cal G})$ of ${\cal G}$ and all edge of $H$ are edge in ${\cal G}.$\\
For $v_1,\,v_n \,\in \,V({\cal G})$ we define a $v_1-v_n$ walk on
${\cal G}$ as a sequence of vertices and edges \linebreak
$W=(v_1,\,e_1,\,v_2,\,e_2,...,v_n)$ such that
$$e_i={\{v_i,\,v_{i+1}\}}\in E(\cal{G})$$ if
$v_1=v_n$ we say that $W$ is a closed walk on ${\cal G},$ this closed walk is said to be a cycle when $e_i \neq e_j,\,\forall\,\, i,j \,\,\in \{1,...,n-1\}.$\\
{\bf degree of a vertex.}\\ Let ${\cal G}$ be a graph and $v \in
V({\cal G})$, the degree of $v$, noted by $p(v)$ is defined as :
$$p(v)=\sharp \{e \in E({\cal G})\,:\, e=\{v,\,w\},\,w \in V\}+\sharp \{e \in E({\cal G})\,:\,e=\{v\} \}$$
{\bf connected graph.}\\A graph ${\cal G}$ is called connected if
there is a $v-w$ walk for all $v,\,w \in V({\cal G}),$ otherwise
${\cal G}$ is
disconnected.\\
{\bf rooted tree.}\\a rooted tree is a pair $(T,\,v)$ such that $T$
is a tree and $v$ a vertex of $T$ with $p(v)=1.$ $v$ is called the root of the tree $T$. \\
{\bf leaf.}\\Let $(T,\,v)$ be a rooted tree, then any vertex $w$ of
$T$ such that $p(w)=1,\,w\neq v,$ is called a leaf of $T.$\\We note
$E(T),\,V(T)$ and $L(T)$ respectively the sets of edges, vertices
and leaves of the rooted tree $T.$\\
{\bf rooted subtree.}\\ A rooted subtree of the rooted tree
$(T,\,x)$ is a pair $(S,\,v)$ where $S$ is a subgraph of $T$ and
$v\in V(T)$ such that $\forall w \in V(S)$ we have $v\leq w$ and $p(v)=1.$\\
A {\bf cut}, ${\mathcal{ C}(T)}$ of the rooted tree $T$ in the vertex $v$,
directly connected to the root by an edge $e=\{x,\,v\}$, is the
uniquely defined collections $\{T_1,...,T_n\}$ of rooted subtrees of
$T$ with root $v$ such that $\bi_{j=1}^n L(T_j)=L(T).$ \\{\bf
Attachment of rooted trees.}\\ Let $x \in \Gamma$ and
$(T_1,\,y)...(T_n,\,y)$ be a collections of rooted trees with
$V(T_i)\,$ and $E(T_i)$ design respectively the set of vertices and
edges of the $i-$th rooted tree. We define an attached rooted tree
$(T,\,x)$ by :\\ $V(T)=\bi_{i=1}^n V(T_i)\cup \{x\}$ and
$E(T)=\bi_{i=1}^n E(T_i) \cup
{\{x,\,y\}}$ respectively the set of vertices and edges of $T$.\\
{\bf Rooted trees with two types of leaves.}\\ Let $(T,\,x)$ be a
rooted tree and $L(T)$ be the set of the leaves of $T$, $(T,\,x)$ is
said to be a rooted tree with two types of leaves if and only if
:\\$L(T)=L_1(T)\cup L_2(T)$ and $l\,\in \,L(T)\setminus L_1(T)
\Longrightarrow l\,\in\,L_2(T).$\\We said that $L_1(T)$ is the set
of the leaves of type one of $T$ and $L_2(T)$ is the set of the
leaves of type two of $T$.\\
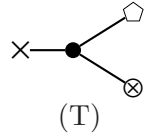
\centerline{\begin{pspicture}(-1,-1)(1,1)
\cnode*(1.2,0){3pt}{B_1}
\dotnode[dotstyle=pentagon,dotscale=2](2,0.5){G'}
\dotnode[dotstyle=otimes,dotscale=2](2,-0.5){G'''}
\dotnode[dotstyle=x,dotscale=2](0.5,0){H_1}\ncline{B_1}{H_1}\ncline{B_1}{G'}\ncline{B_1}{G''}\ncline{B_1}{G'''}
\end{pspicture}(T)}	
\\
In this graph (T) the inner vertices are {\bf distinguishable} and
have {\bf non distinguishable} legs :\\Let $x$ be the root of the
tree $T$ and we consider the edges $e_1=\{y_1,\,z_1\}$ and
$e_2=\{y_1,\,z_2\}$, the value of the tree $T$ [cf. Def.(3.1)], is
the same when we permute the edges i.e when $e_1=\{y_1,\,z_2\}$ and
$e_2=\{y_1,\,z_1\}.$\\
In this Parisi-Wu graph $({\cal{G}})$,with one inner vertex $y_1$
and two roots $x_1,\,x_2,$ the empty vertex have {\bf non
distinguishable} legs .\\Let $z$ be the empty vertex and
$e_1=\{y_1,\,z\}$ , $e_2=\{y_2,\,z\}$ be the edges connected to the
empty vertex.\\The value of the generalized Parisi-Wu graph
$({\cal{G}})$, [cf. Def.(4.4)] is the same when we permute the edges
i.e when $e_1=\{y_2,\,z\}$ and $e_2=\{y_1,\,z\}.$ \vspace{1cm}

\noindent{\bf Acknowledgements:} The authors thanks
Professor Sergio Albeverio for his interest and for valuable discussions. 

\vspace{1cm}
\noindent {\sc Hanno Gottschalk and Boubaker Smii \\
\rm Institut f\"ur angewandte Mathematik\\
Wegeler str. 6\\
D-53115 Bonn, Germany\\
e-mails: gottscha/boubaker@wiener.iam.uni-bonn.de }
\end{document}